\documentclass[a4paper,11pt]{article}
\usepackage[top=30truemm,bottom=30truemm,left=30truemm,right=30truemm]{geometry}
\usepackage{amsmath}
\usepackage{amssymb}
\usepackage{float}
\usepackage{color}
\usepackage[dvipdfmx]{graphicx}
\usepackage{graphics}
\usepackage{graphicx}
\usepackage{tcolorbox}
\usepackage{latexsym}
\usepackage{ascmac}
\usepackage{framed}
\usepackage{fancybox}
\usepackage{theorem}
\usepackage{stmaryrd}
\usepackage{here}
\usepackage{mathtools}

\mathtoolsset{showonlyrefs=true} 

\theoremstyle{break}
\def\qed{\hfill$\Box$}

\newtheorem{defn}{Definition}

\newtheorem{prop}{Proposition}
\newtheorem{thm}{Theorem}
\newtheorem{cor}{Corollary}

\newtheorem{rem}{Remark}                                   

\makeatletter
\@addtoreset{equation}{section}

\makeatother

\makeatletter
\@addtoreset{figure}{section}
\def\thefigure{\thesection.\arabic{figure}}
\makeatother

\title{
Classification of nonnegative traveling wave solutions for certain 1D degenerate parabolic equation and porous medium equation 
}

\author{
Yu Ichida
\thanks{Meiji University Organization for Strategic Coordination of Research and Intellectual Properties, 1-1 Kanda-Surugadai, Chiyoda-ku, Tokyo 101-8301, Japan, {\tt ichidayu@meiji.ac.jp} } $^{,}$ \footnote{JSPS Research Fellow (Research Fellow of Japan Society for the Promotion of Science)}
\,\, , 
Takashi Okuda Sakamoto
\thanks{Graduate School of Science and Technology, Meiji University, 1-1-1, Higashimita Tama-ku Kawasaki Kanagawa 214-8571, Japan, {\tt sakamoto@meiji.ac.jp}}
}

\begin{document}
\maketitle
\begin{abstract}
This paper reports results on the classification of traveling wave solutions, including nonnegative weak sense, in the spatial 1D degenerate parabolic equation.
These are obtained through dynamical systems theory and geometric approaches (in particular, Poincar\'e compactification).
Classification of traveling wave solutions means enumerating those that exist and presenting properties of each solution, such as its profile and asymptotic behavior.
The results examine a different range of parameters included in the equation, using the same techniques as discussed in the earlier work [Y. Ichida, Discrete Contin. Dyn. Syst., Ser. B, {\bf{28}} (2023), no. 2, 1116--1132].
In a clear departure from this previous work, the classification results obtained in this paper and the successful application of known transformation also yield results for the classification of (weak) nonnegative traveling wave solutions for spatial 1D porous medium equations with special nonlinear terms and the simplest porous medium equation.
Finally, the bifurcations at infinity occur in the two-dimensional ordinary differential equations that characterize these traveling wave solutions are shown.
\end{abstract}


\begin{center}
\end{center}

\section{Introduction}
\label{sec:pDNPE1}
In this paper, we consider the classification of nonnegative (weak sense) traveling wave solutions of the following three parabolic partial differential equations.
The first is the spatial one-dimensional degenerate parabolic equation:
\begin{equation}
U_{t}=U^{p}(U_{xx}+\mu U)-\delta U, \quad t>0, \quad x\in \mathbb{R}
\label{eq:pDNPE1-1}
\end{equation}
with $0<p<1$, $\mu>0$ and $\delta=0$ or $1$.
The second is the spatial one-dimensional porous medium equation
\begin{equation}
V_{t}=(V^{m})_{xx}+\mu V^{m} -\delta (1-p)V, \quad t>0, \quad x\in \mathbb{R}, \quad \quad m=\dfrac{1}{1-p}>1
\label{eq:pDNPE1-5}
\end{equation} 
with the special linear and nonlinear terms.
\eqref{eq:pDNPE1-5} is the equation obtained by applying the transformation
\begin{equation}
U(t, x) = (1-p)^{-\frac{1}{p}} (V(t,x))^{\frac{1}{1-p}}
\label{eq:pDNPE1-4}
\end{equation}
to \eqref{eq:pDNPE1-1}, which is also introduced in Winkler \cite{Wink03}.
Finally, the third is the simple spatial one-dimensional porous medium equation
\begin{equation}
V_{t}=(V^{m})_{xx}, \quad t>0, \quad x\in \mathbb{R}, \quad m>1.
\label{eq:pDNPE1-po}
\end{equation}
This is obtained by $\delta=0$ and $\mu=0$ in \eqref{eq:pDNPE1-5}.

First, a brief background of equation \eqref{eq:pDNPE1-1} is given.
\eqref{eq:pDNPE1-1} for $\delta=0$ comes from the time evolution problem of a certain plane curve.
Bordering on $p=1$, the case $0<p<1$ is related to the problem of the movement of an expanding curve, while the case $p>1$ is related to the problem of curve shortening (see, for instance, \cite{LPT12} and references therein).
According to \cite{Anada, Poon} and references therein, \eqref{eq:pDNPE1-1} for $\delta=0$ arises in the modeling of the solar flares in astrophysics and the resistive diffusion of a force-free magnetic field in a plasma confined between two walls.
Although \eqref{eq:pDNPE1-1} with  $\delta=1$ has no background on specific phenomena, it plays an important role in understanding the solution structure of the equation, as described next.

Next, the mathematical problem in \eqref{eq:pDNPE1-1} is explained.
There have been many studies on blow-up solutions in \eqref{eq:pDNPE1-1} for the case $\delta=0$ (see \cite{Anada, AIU22, ange91, Ange, LPT12, Poon, Wink03}).
According to these references, the solution has the blow-up of Type I when $0<p<2$.
When $p\ge2$, the solution has the blow-up of Type II.
In particular, much work has been done in recent years on the derivation of blow-up rates in $p\ge 2$ (see, for instance, \cite{Anada, Poon, AIU22} and references therein).
In these studies, the following equation, which restricts the range of $(t,x)$ considered in \eqref{eq:pDNPE1-1} with $\delta=0$, is used:
\begin{equation}
\bar{U}_{t}=\bar{U}^{p}( \bar{U}_{xx}+ \mu\bar{U}), \quad (t,x)\in (0,T)\times (-L,L), 
\label{eq:Known1} 
\end{equation}
where $p \in \mathbb{R}$ and $T<\infty$, and impose Dirichlet or periodic boundary conditions.
Note that the unknown function is $\bar{U}=\bar{U}(t,x)$ to distinguish it from \eqref{eq:pDNPE1-1}.
\cite{AIU22, Ange, DNPE, Poon, Wink03}, $\tilde{U}(\tau,x)$ satisfies
\begin{equation}
\tilde{U}_{\tau}=\tilde{U}^{p}(v_{xx}+\mu \tilde{U}-\tilde{U}^{-p+1}), \quad (\tau,x)\in (0,+\infty)\times (-L,L)
\label{eq:Known2} 
\end{equation}
by introducing the following rescaled function:
\begin{equation}
\tilde{U}(\tau,x)=(pT)^{\frac{1}{p}}e^{-\tau} \bar{U}(t,x),\quad t=T(1-e^{-p\tau}),\quad \tau\in (0,+\infty).
\label{eq:pDNPEt1}
\end{equation}
This transformation is derived from the self-similar solution and moves the blow-up time to infinity.
Note that in this paper, \eqref{eq:Known1} and \eqref{eq:Known2} are combined and \eqref{eq:pDNPE1-1} is considered as the whole space with respect to space.
The reason for considering the whole space is related to the fact that the blow-up solution is not considered in this study and that the proof of the main results, which will be discussed later, is obtained by focusing only on the structure of the equations.
In \cite{Ange, Poon, AIU22}, the behavior of the traveling wave solution of \eqref{eq:Known2} is investigated to evaluate $\tilde{U}(\tau, 0)$, which plays an important role in studying the blow-up rate.
It can be seen that these previous studies have successfully derived a lower bound for the blow-up rate by using the information from the traveling wave solution in \eqref{eq:Known2}.
According to \cite{Ange}, the traveling wave solution in \eqref{eq:Known2} corresponds to examining a special self-similar solution of \eqref{eq:Known1}.
More precisely, the traveling wave solutions correspond to special self-similar solutions of curve shortening which evolves by rotating and contracting simultaneously.

Next, we briefly review the results for the traveling wave solution of \eqref{eq:Known2} and  \eqref{eq:pDNPE1-1}.
Let $\tilde{U}(\tau, x)=U_{c}(\xi)$ ($\xi=x-ct$) be the traveling wave solution of \eqref{eq:Known2}.
Angenent \cite{ange91} and Angenent-Vel\'azquez \cite{Ange} showed that in \eqref{eq:Known2} for $p=2$ and $\mu=1$, there exists a unique traveling wave solution $U_{c}(\xi)$ such that the property:
\[
U_{c}'(0)=0, \quad U_{c}'(\xi)>0 \quad {\rm{for}} \quad \xi<0, \quad U_{c}(\xi)\to 0 \quad {\rm{as}} \quad \xi \to -\infty
\]
exists.
In Lin-Poon-Tsai \cite{LPT12} and Poon \cite{Poon}, this result is extended to the case $p>2$ and $\mu=1$.
These are shown in the phase plane analysis method.
Motivated by these studies, Ichida-Sakamoto \cite{DNPE} considers \eqref{eq:pDNPE1-1} in the whole of space 1D and gives a classification of traveling wave solutions, including weak meaning for the special $p$.
More specifically, it gives information about the existence, shape, and asymptotic behavior of the (weak) traveling wave solution of \eqref{eq:pDNPE1-1} in $p \in 2\mathbb{N}$, $\mu=1$, and both $\delta=0$ and $\delta=1$.
Note that this result encompasses the results of the previous studies \cite{ange91, Ange, LPT12, Poon} mentioned above.
In \cite{DNPE}, the main results are given by applying dynamical systems theory, Poincar\'e compactification (see below for details), and a geometric technique for the desingularization of vector fields called the blow-up technique (see \cite{AFJ, MB, FAL}).
The meaning of the classification of solutions derives from the fact that the above method reveals all dynamics of the two-dimensional system of ordinary differential equations (for short, ODEs) obtained by the traveling wave coordinates, including to infinity.
Furthermore, it should be noted that the weak traveling wave solution obtained for $\delta=1$ in Theorem 2 of \cite{DNPE} is partially utilized in the discussion of \cite{AIU22}.
As mentioned in \cite{cDNPE}, Ichida-Matsue-Sakamoto \cite{DNLA} gave a refined asymptotic behavior, which was not obtained in the preceding work \cite{DNPE}, by an appropriate asymptotic study and properties of the Lambert $W$ function.
In \cite{DNPE, DNLA}, it is necessary to assume $1<p \in \mathbb{N}$ from previous studies \cite{AFJ, MB, FAL}, since the discussion process uses the method of blow-up technique.
Therefore, there is no discussion on the general $1<p \in \mathbb{R}$.
Therefore, the author \cite{cDNPE} investigates the nonnegative (weak) traveling wave solution of the equation
\begin{equation}
u_{t}=u u_{xx}-\gamma(u_{x})^{2}+ku^{2}-\delta p u, \quad t>0, \quad x\in \mathbb{R},
\label{eq:pDNPE1-3}
\end{equation}
for $u$ obtained by introducing the transformation
\begin{equation}
u(t,x) = (U(t,x))^{p}, \quad \gamma=\dfrac{p-1}{p}, \quad k=p\mu. 
\label{eq:pDNPE1-2}
\end{equation}
Then, the classification of traveling wave solutions of \eqref{eq:pDNPE1-3} can be obtained by the same argument as for \cite{DNPE}.
For this result, by using $U=u^{1/p}$, which is the inverse transformation of \eqref{eq:pDNPE1-2}, we can obtain the classification of nonnegative traveling wave solutions in both $\delta=0$ and $\delta=1$ for $1<p \in \mathbb{R}$ in \eqref{eq:pDNPE1-1}. 
The key to this argument is that in the classification of traveling wave solutions of \eqref{eq:pDNPE1-3} obtained by \eqref{eq:pDNPE1-2}, it is not necessary to use the blow-up technique since $p$ appears in the coefficient part instead of the exponent in the calculation process.
See \cite{cDNPE} for more details.
This result generalizes the \cite{DNPE, DNLA} result to $1<p \in \mathbb{R}$.

The methods (in particular, the Poincar\'e compactification) employed in the above previous studies \cite{DNPE, cDNPE} are briefly described.
The Poincar\'e compactification is the key method in this paper.
See \cite{FAL, QTW, DNPE, cDNPE, Matsue1, Matsue2} for details and geometric images.
This is one of the compactifications of the original phase space, the embedding of $\mathbb{R}^{n}$ into $\mathbb{R}^{n+1}$ in the unit's upper hemisphere.
The procedure of this method is briefly described.
As described in \cite{Matsue1}, this compactification makes infinity in the original phase space correspond to the boundary of the compact manifold.
The infinity is divided into several parts, each of which is projected to the local coordinate corresponding to infinity.
The dynamics at each local coordinate are then examined, i.e., the dynamics including a part of the segmented infinity.
By combining the information obtained in these local coordinates, it is possible to obtain the dynamics including infinity in the original phase space (hereafter referred to as the dynamics on the Poincar\'e disk).
This method has been used, for instance, in the analysis of the Li\'enard equation (\cite{FAL} and references therein) and in the reconstruction of blow-up solutions of ODEs in the view of dynamical systems theory (see \cite{Matsue1, Matsue2}).

Originally, the traveling wave solution is considered in terms of its contribution to the clarification of the blow-up rate.
It is pointed out that in terms of dynamical systems theory and geometric approaches (especially Poincar\'e compactification) in the discussion of \cite{DNPE, DNLA, cDNPE}, it is possible to obtain a rich property of previously unknown information on traveling wave solutions.
In this paper, we attempt to classify nonnegative traveling wave solutions in \eqref{eq:pDNPE1-1} for $0<p<1$, which has not been clarified before.
Note that this is a problem of purely mathematical interest, not an analysis of traveling wave solution  for more accurate blow-up rate derivation.
Even if \eqref{eq:pDNPE1-1} is directly attributed to the problem of investigating the behavior of two-dimensional ODEs by introducing the traveling wave coordinates in the same as \cite{DNPE, DNLA}, it is a nontrivial problem that cannot be investigated for $0<p<1$ in an analysis such as \cite{DNPE, DNLA}.
Therefore, we first consider the problem of classifying nonnegative traveling wave solutions of \eqref{eq:pDNPE1-3} in $0<p<1$ using the same transformation \eqref{eq:pDNPE1-2} also introduced in \cite{cDNPE}.
Then, as a similar argument to \cite{cDNPE}, we investigate the dynamics of a system of two-dimensional ODEs including to infinity satisfied by the traveling wave coordinates by Poincar\'e compactification.
By using $U=u^{1/p}$, the result in the equation $U$ is obtained.
Although similar to the problem of \cite{cDNPE}, it should be emphasized that the work in this paper is significant in that nonnegative traveling wave solutions in nontrivial $0<p<1$ can be investigated in the same way as $1<p \in \mathbb{R}$ by introducing the transformation \eqref{eq:pDNPE1-2}.
For the method, we use the Poincar\'e compactification used in \cite{DNPE, DNLA, cDNPE}.

Comparing \cite{cDNPE} with this paper means that Poincar\'e compactification and dynamical systems theory methods are effective in investigating traveling wave solutions even if the range of $p$ in \eqref{eq:pDNPE1-1} is different.
It should be emphasized here that the classification of traveling wave solutions for \eqref{eq:pDNPE1-1} obtained in this paper leads to applications that are clearly different from those of \cite{cDNPE}.
As mentioned at the beginning, when $0<p<1$, we apply \eqref{eq:pDNPE1-4}, which is also introduced in Winkler \cite{Wink03}, we then obtain \eqref{eq:pDNPE1-5}.
\eqref{eq:pDNPE1-5} is a kind of porous medium equation.
As described in \cite{XuYin}, finite propagation is known as a fundamental feature of degenerate diffusion equations.
Therefore, it is one of the issues to clarify the existence and various properties of traveling wave solutions, which are characteristic solutions.
See, for instance, \cite{PaVa} for a study of traveling wave solutions of \eqref{eq:pDNPE1-5} with $\delta=0$.
In \cite{PaVa}, results on the global existence and finite propagation of traveling wave solutions are also presented.
However, no results on the classification of nonnegative traveling wave solutions have been given.
In this paper, as an application of the results on the classification of nonnegative (weak) traveling wave solutions in \eqref{eq:pDNPE1-1} for $0<p<1$, we obtain the classification of nonnegative (weak) traveling wave solutions in \eqref{eq:pDNPE1-5}.
In \eqref{eq:pDNPE1-5}, by setting $\mu=0$ and $\delta=0$, we obtain \eqref{eq:pDNPE1-po} and the classification of nonnegative traveling wave solutions of the simplest porous medium equation.
For background on the simplest porous medium equation and its mathematical problems, see, for instance, \cite{Aron, XuYin} and references therein.

For the two-dimensional ODEs obtained by introducing the traveling wave coordinate into \eqref{eq:pDNPE1-po}, it is difficult to correctly extract the structure at infinity by applying Poincar\'e compactification.
In other words, it is difficult to obtain a classification of traveling wave solutions of the space one-dimensional simple porous media equation with the same arguments as for \cite{QTW, DNPE, cDNPE}.
In this sense, it is a meaningful result (Corollary \ref{cor:pDNPE2-3}) as a system derived from Theorem \ref{thm:pDNPE2-1} in this paper.

Next, we derive the two-dimensional ODE system guided by the traveling wave coordinates that play a central role in the discussion before discussing the main results.
Then, for \eqref{eq:pDNPE1-3}, we introduce the following change of variables:
\[
\phi(\xi)=u(t,x), \quad \xi=x-ct, \quad 0<c\in \mathbb{R}.
\]
The equation of $\phi(\xi)$ is then reduced to 
\begin{equation}
\phi \phi''=-c\phi'+\gamma (\phi')^{2}-k\phi^{2}+\delta p\phi, 
\quad \left( \, '=\dfrac{d}{d\xi}, \quad ''=\dfrac{d^{2}}{d\xi^{2}}  \,\right).
\label{eq:pDNPE1-6}
\end{equation}
Then, \eqref{eq:pDNPE1-6} is equivalent to 
\begin{equation}
\begin{cases}
\phi'=\psi,
\\
\psi'=-c\phi^{-1}\psi+\gamma \phi^{-1}\psi^{2}-k\phi+\delta p,
\end{cases}
\quad \left( \, '=\dfrac{d}{d\xi}\, \right).
\label{eq:pDNPE1-7}
\end{equation}
This is also derived in \cite{cDNPE}, which is the same in that $k>0$, but differs in that $\gamma<0$.
Similarly to \cite{cDNPE}, all dynamics of \eqref{eq:pDNPE1-7} up to infinity can be obtained by Poincar\'e compactification.
Then, all connecting orbits of \eqref{eq:pDNPE1-7} including to infinity can be classified.
This result indicates that the nonnegative traveling wave solutions corresponding to each connecting orbit in the Poincar\'e disk are classified.
Therefore, once the existence of each traveling wave solution is known, its shape and asymptotic behavior can be obtained.

By organizing the results of \cite{cDNPE} and this paper, we can also investigate the change in stability of the finite equilibria and equilibria at infinity of \eqref{eq:pDNPE1-7} for $p>0$.
See Definition 3.13 in \cite{Matsue1} for the definition of equilibria at infinity.
In Section \ref{sec:pDNPE5}, if we move $p$ from $p>1$ to $0<p<1$ after $p=1$, we can observe the bifurcations of the equilibria at infinity.
This is obtained by observing the equilibria on the $\psi$-axis up to infinity in Poincar\'e disk and the change in $p$.
Since \eqref{eq:pDNPE1-7} characterizes the traveling wave solutions of \eqref{eq:pDNPE1-3} and \eqref{eq:pDNPE1-1}, therefore, the bifurcations of equilibria at infinity in ODE is related to the qualitative change of the non-negative traveling wave solutions (cf. \cite{cDNPE}).

The paper is organized as follows.
The next section describes the main results of this paper.
In Section \ref{sec:pDNPE3}, the dynamics of \eqref{eq:pDNPE1-7} on the Poincar\'e disk are obtained based on Poincar\'e compactification and classical dynamical systems theory.
The discussion is almost the same as for \cite{cDNPE}, but for the reader's convenience, we briefly discuss its details and differences from \cite{cDNPE}.
The proof of the main result is completed in Section \ref{sec:pDNPE4}.
Finally, in Section \ref{sec:pDNPE5}, we summarize the results of this paper and those of \cite{cDNPE} to describe changes in the stability of the equilibria  including infinity of the two-dimensional ODEs, i.e., bifurcaition that characterize the traveling wave solutions.

\section{Main results}
\label{sec:pDNPE2}
Before describing the main results in this paper, we define the necessary terms, concepts, and symbols.
These definitions are not general terms but have already been adopted by \cite{QTW, DNPE, cDNPE}.
In this paper, we adopt these definitions in order to compare the obtained results with those of \cite{cDNPE}.

\begin{defn}
\label{def:pDNPE2-1}
We say that a function $u(t,x) \equiv \phi(\xi)$ is a quasi traveling wave of \eqref{eq:pDNPE1-3} if the function $\phi(\xi)$ is a solution of \eqref{eq:pDNPE1-6} on a finite interval or semi-infinite
interval.
It is defined similarly in \eqref{eq:pDNPE1-1}, \eqref{eq:pDNPE1-5}, and \eqref{eq:pDNPE1-po}.
\end{defn}

\begin{defn}
\label{def:pDNPE2-2}
We say that a function $u(t,x) \equiv \phi(\xi)$ is a quasi traveling wave with quenching of \eqref{eq:pDNPE1-3}
if the function $u(t,x)$ is a quasi traveling wave of \eqref{eq:pDNPE1-3} on a finite (resp. semi-infinite interval) such that $u=\phi$ reaches $0$ and $|u_{\xi}|=|\phi'|=|\psi|$ becomes infinity at finite endpoints (resp. point) of the finite (resp. semi-infinite) interval. 
More precisely, the function $\phi(\xi)$ is a solution of \eqref{eq:pDNPE1-6} on a finite (resp. semi-infinite) interval $(\xi_{-}, \xi_{+})$ (resp.$( \xi_{-},+\infty)$) and $\phi(\xi)$ satisfies $\phi(\xi) \in C^{2}(\xi_{-},\xi_{+}) \cap C^{0}[\xi_{-},\xi_{+}]$, $-\infty<\xi_{-}<\xi_{+}<+\infty$ (resp. $\phi(\xi) \in C^{2}(\xi_{-},+\infty) \cap C^{0}[\xi_{-},+\infty)$, $|\xi_{-}|<\infty$).
In addition, the function satisfies 
\begin{align}
&
\lim_{\xi \nearrow \xi_{+} - 0} \phi(\xi) = 0, \quad 
\lim_{\xi \nearrow \xi_{+} - 0} |\phi'(\xi)| = \infty,
\label{eq:pDNPE2-1}
\\ 
&
\lim_{\xi \searrow \xi_{-} + 0} \phi(\xi) = 0, \quad
\lim_{\xi \searrow \xi_{-} +0} |\phi'(\xi)| = \infty
\label{eq:pDNPE2-2}
\end{align}
(resp. only \eqref{eq:pDNPE2-2}).
It is defined similarly in \eqref{eq:pDNPE1-1}, \eqref{eq:pDNPE1-5}, and \eqref{eq:pDNPE1-po}.
\end{defn}

\begin{defn}
\label{def:pDNPE2-3}
We say that a function $u(t,x) \equiv \phi(\xi)$ is a quasi traveling wave with the singularity of \eqref{eq:pDNPE1-3} 
if the function $u(t,x)$ is a quasi traveling wave of \eqref{eq:pDNPE1-3} on a semi-infinite interval such that $\phi$ reaches $0$ and only left differentiation is possible and it becomes a constant at the finite end point of the semi-infinite interval. 
More precisely, the function $\phi(\xi)$ is a solution of \eqref{eq:pDNPE1-6} on a semi-infinite interval $(-\infty, \xi_{+})$ ($\phi(\xi) \in C^{2}(-\infty, \xi_{+})$, $|\xi_{+}|<\infty$), and satisfies 
\begin{equation}
\lim_{\xi \nearrow \xi_{+} -0} \phi(\xi) = 0
\quad {\rm and} \quad
\lim_{\xi \nearrow \xi_{+} -0} \phi'(\xi) = -C 
\label{eq:pDNPE2-3}
\end{equation}
with $C>0$.
This definition holds for the semi-infinite interval $(\xi_{-}, +\infty)$ and finite interval $(\xi_{-}, \xi_{+})$ as well.
\end{defn}

A function that satisfies \eqref{eq:pDNPE2-2} at one endpoint of a finite interval and \eqref{eq:pDNPE2-3} at the other endpoint is newly called a quasi traveling wave with quenching and singularity.
Note that this concept does not appear in \cite{DNPE, cDNPE}, but is a combination of those defined in \cite{DNPE, cDNPE} and is the first term defined in this paper.

\begin{defn}
\label{def:pDNPE2-4}
Let $u(\xi)$ be a quasi traveling waves with quenching of \eqref{eq:pDNPE1-3} on a finite interval.
Then, we say that a function 
\[
u^{*}(\xi) = 
\begin{cases}
u(\xi), & \xi \in (\xi_{-}, \xi_{+}) ,
\\
0, & else
\end{cases}
\]
is a weak traveling wave solution (with quenching) of \eqref{eq:pDNPE1-3}.
If $u(\xi)$ is a quasi traveling wave with the singularity (resp. quasi traveling wave with quenching and singularity) \eqref{eq:pDNPE1-3} on a finite or semi-infinite interval, then a weak traveling wave solution (with the singularity (resp. quenching and singularity)) of \eqref{eq:pDNPE1-3} is defined in the same way.
They are defined similarly in \eqref{eq:pDNPE1-1}, \eqref{eq:pDNPE1-5}, and \eqref{eq:pDNPE1-po}.
\end{defn}

The above definition implies that $u^{*}(\xi)$ satisfies
\[
\int_{\mathbb{R}} \left[u\varphi_{\xi}(c-u_{\xi})+(\gamma+1)(u_{\xi})^{2}\varphi-u^{2}(\varphi_{\xi\xi}+k\varphi)+\delta pu\varphi \right]\, d\xi=0
\]
for all $\varphi\in C_{0}^{\infty}(\mathbb{R})$.
Note that this is the same as \cite{cDNPE}.

Under these definitions, the main results of this paper are presented.
Note that $\phi(\xi)=u(\xi)$ and $\phi'(\xi)=\psi(\xi)=u_{\xi}(\xi)$ hold.
Hereinafter, note that the meaning of the symbol  $f(\xi)\sim g(\xi)$ as $\xi\to a$ is as follows:
\[
\lim_{\xi\to a}\left| \dfrac{f(\xi)}{g(\xi)} \right|=1.
\]

First, we present results on the classification of nonnegative traveling wave solutions of \eqref{eq:pDNPE1-3} (including the weak sense) for $\delta=0$ and $\delta=1$, respectively.

\begin{prop}
\label{prop:pDNPE2-1}
Assume that $0<p<1$, $\gamma<0$, $k>0$, $\mu>0$, and $\delta=0$.
Then, for a given positive constant $c$, the equation \eqref{eq:pDNPE1-3} has following three type of weak traveling wave solutions:

\begin{enumerate}
\item[(I)] 
There exists a family of weak traveling wave solutions (with quenching) such that it corresponds to the family of orbits in \eqref{eq:pDNPE1-7}.
Each solution $u(\xi)$ satisfies the following:
\begin{enumerate}
\item[(I1)]
$\displaystyle \lim_{\xi \searrow \xi_{-} + 0} u(\xi) = \lim_{\xi \to +\infty} u(\xi) = \lim_{\xi \to +\infty} u'(\xi) =  0$,
$\displaystyle\lim_{\xi \searrow \xi_{-} + 0} u'(\xi) = +\infty$.
\item[(I2)]
$u(\xi)>0$ holds for $\xi\in (\xi_{-},+\infty)$ and $u(\xi)=0$ holds for $\xi\in (-\infty, \xi_{-}]$.
\item[(I3)]
There exists a constant $\xi_{0}\in (\xi_{-}, +\infty)$ such that the following holds:
$u'(\xi)>0$ for $\xi\in (\xi_{-}, \xi_{0})$, $u'(\xi_{0})=0$ and $u'(\xi)<0$ for $\xi\in (\xi_{0}, +\infty)$.
\end{enumerate}
In addition, the asymptotic behavior of $u(\xi)$ and $u'(\xi)$ for $\xi \searrow \xi_{-}+0$ are
\begin{equation}
\begin{cases}
u(\xi) \sim A_{1}(\xi-\xi_{-})^{p},
\\
u'(\xi) \sim A_{2}(\xi-\xi_{-})^{p-1}
\end{cases}
\quad {\rm{as}} \quad \xi \searrow \xi_{-}+0,
\label{eq:pDNPE2-4}
\end{equation}
where $A_{j}$ are positive constants, and the asymptotic behavior of $u(\xi)$ and $u'(\xi)$ for $\xi \to +\infty$ are
\begin{equation}
\begin{cases}
u(\xi) \sim ck^{-1}\xi^{-1},
\\
u'(\xi) \sim -ck^{-1}\xi^{-2},
\end{cases}
\quad {\rm{as}} \quad \xi\to +\infty.
\label{eq:pDNPE2-5}
\end{equation}

\item[(II)] 
There exists a family of weak traveling wave solutions (with quenching) such that it corresponds to the family of orbits in \eqref{eq:pDNPE1-7}.
Each solution $u(\xi)$ satisfies the following:
\begin{enumerate}
\item[(II1)]
$\displaystyle \lim_{\xi \searrow \xi_{-} + 0} u(\xi) = \lim_{\xi \nearrow \xi_{+}-0} u(\xi) =0$,
$\displaystyle \lim_{\xi \to \xi_{-}+0}  u'(\xi) =  +\infty$,
$\displaystyle\lim_{\xi \nearrow \xi_{+} - 0} u'(\xi) = -\infty$.
\item[(II2)]
$u(\xi)>0$ holds for $\xi\in (\xi_{-},\xi_{+})$ and $u(\xi)=0$ holds for $\xi\in (-\infty, \xi_{-}] \cup [\xi_{+}, +\infty)$.
\item[(II3)]
There exists a constant $\xi_{0}\in (\xi_{-}, \xi_{+})$ such that the following holds:
$u'(\xi)>0$ for $\xi\in (\xi_{-}, \xi_{0})$, $u'(\xi_{0})=0$ and $u'(\xi)<0$ for $\xi\in (\xi_{0}, \xi_{+})$.
\end{enumerate}
In addition, the asymptotic behavior of $u(\xi)$ and $u'(\xi)$ for $\xi \searrow \xi_{-}+0$ are expressed as \eqref{eq:pDNPE2-4}, and the asymptotic behavior of $u(\xi)$ and $u'(\xi)$ for $\xi \nearrow \xi_{+}-0$ are
\begin{equation}
\begin{cases}
u(\xi) \sim A_{3}(\xi_{+}-\xi)^{p},
\\
u'(\xi) \sim -A_{4}(\xi_{+}-\xi)^{p-1}
\end{cases}
\quad {\rm{as}} \quad \xi \nearrow \xi_{+}-0,
\label{eq:pDNPE2-6}
\end{equation}
where $A_{j}$ are positive constants.

\item[(III)] 
There exists a weak traveling wave solution (with quenching and singularity) such that it corresponds to the orbit in \eqref{eq:pDNPE1-7}.
The solution $u(\xi)$ satisfies (II2), (II3), and the following:
\begin{enumerate}
\item[(III1)]
$\displaystyle \lim_{\xi \searrow \xi_{-} + 0} u(\xi) = \lim_{\xi \nearrow \xi_{+}-0} u(\xi) =0$,
$\displaystyle \lim_{\xi \searrow \xi_{-}+0}  u'(\xi) =  +\infty$,
$\displaystyle\lim_{\xi \nearrow \xi_{+} - 0} u'(\xi) = -C$ with $C>0$. 
\end{enumerate}
In addition, the asymptotic behavior of $u(\xi)$ and $u'(\xi)$ for $\xi \searrow \xi_{-}+0$ are expressed as \eqref{eq:pDNPE2-4}, and the asymptotic behavior of $u(\xi)$ and $u'(\xi)$ for $\xi \nearrow \xi_{+}-0$ are
\begin{equation}
u(\xi) \sim A_{1}(\xi_{+}-\xi),
\quad
u'(\xi) \sim -C
\quad {\rm{as}} \quad \xi \nearrow \xi_{+}-0
\label{eq:pDNPE2-7}
\end{equation}
with a positive constant $C>0$.
\end{enumerate}
\end{prop}

\begin{prop}
\label{prop:pDNPE2-2}
Assume that $0<p<1$, $\gamma<0$, $k>0$, $\mu>0$, and $\delta=1$.
Then, for a given positive constant $c$, the equation \eqref{eq:pDNPE1-3} has three types of weak traveling wave solutions and one traveling wave solution including (II) and (III) in Proposition \ref{prop:pDNPE2-1}.
The remaining two characterizations are as follows:

\begin{enumerate}
\item[(IV)] 
There exists a family of weak traveling wave solutions (with quenching) such that it corresponds to the family of orbits in \eqref{eq:pDNPE1-7}.
Each solution $u(\xi)$ satisfies (I2), (I3), and the following:
\begin{enumerate}
\item[(IV1)]
$\displaystyle 
\lim_{\xi \searrow \xi_{-} + 0} u(\xi)=0$,
$\displaystyle
\lim_{\xi \to +\infty} u(\xi) = \mu^{-1}$,
$\displaystyle
\lim_{\xi \searrow \xi_{-} + 0} u'(\xi) = +\infty$.
\end{enumerate}
In addition, the asymptotic behavior of $u(\xi)$ and $u'(\xi)$ for $\xi \searrow \xi_{-}+0$ are expressed as \eqref{eq:pDNPE2-4}, and the asymptotic behavior of $u(\xi)$ for $\xi \to +\infty$ is
\begin{equation}
u(\xi) \sim \dfrac{1}{\mu} \sim
\begin{cases}
B_{1}e^{\omega_{1}\xi}+B_{2}e^{\omega_{2}\xi}+\dfrac{1}{\mu}, \quad (D>0),
\\
(B_{3}\xi+B_{4})e^{\omega\xi}+\dfrac{1}{\mu},\quad (D=0),
\\
e^{-\frac{\mu c}{2}\xi}\bar{Z}(\xi) +\dfrac{1}{\mu},\quad (D<0),
\end{cases}
\label{eq:pDNPE2-8}
\end{equation}
where $B_{j}$ ($1\le j\le 4$) are constants and
\begin{align*}
&\omega_{1}=\dfrac{-\mu c+\sqrt{D}}{2}<0, \quad \omega_{2}=\dfrac{-\mu c-\sqrt{D}}{2}<0, \quad \omega=-\dfrac{\mu c}{2}<0, \quad D=\mu^{2}c^{2}-4k. 
\\
&\bar{Z}(\xi) = B_{5}\cdot  \sin [\frac{\sqrt{|D|}}{2}\xi]+ B_{6}\cdot \cos[ \frac{\sqrt{|D|}}{2}\xi] 
\end{align*}
with constants $B_{5, 6}$.

\item[(V)] 
There exists a traveling wave solution such that it corresponds to the orbit in \eqref{eq:pDNPE1-7}.
The solution $u(\xi)$ satisfies the following:
\begin{enumerate}
\item[(V1)]
$\displaystyle 
\lim_{\xi \to -\infty} u(\xi) =0,
\lim_{\xi \to +\infty} u(\xi) = \mu^{-1}.
$
\item[(V2)]
$u(\xi)>0$ holds for $\xi\in \mathbb{R}$.
\end{enumerate}
In addition, the asymptotic behavior of $u(\xi)$ for $\xi \to \infty$ is expressed as \eqref{eq:pDNPE2-8}, and the asymptotic behavior of $u(\xi)$ for $\xi \to -\infty$ is
\begin{equation}
u(\xi) \sim \dfrac{Mc^{2}e^{\frac{p}{c}\xi}}{M(\mu c^{2}+1)e^{\frac{p}{c}\xi}-1}
\quad {\rm{as}}\quad \xi\to -\infty,
\label{eq:pDNPE2-9} 
\end{equation}
where $M<0$ is a constant that depends on the initial state $\phi(0)=\phi_{0}$.
\end{enumerate}
\end{prop}

\begin{figure}[t]
\centering
\includegraphics[width=4.6cm]{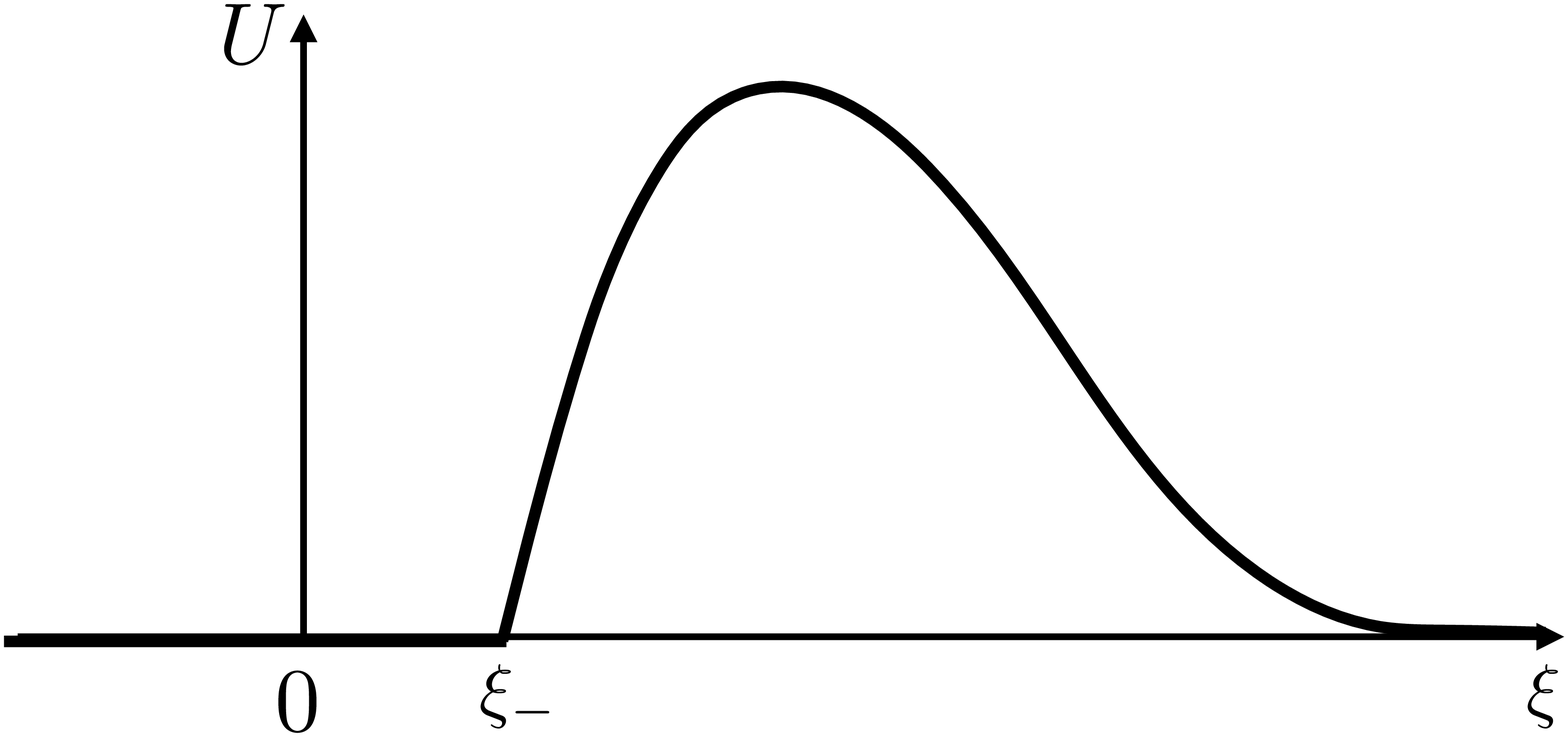}
\includegraphics[width=4.6cm]{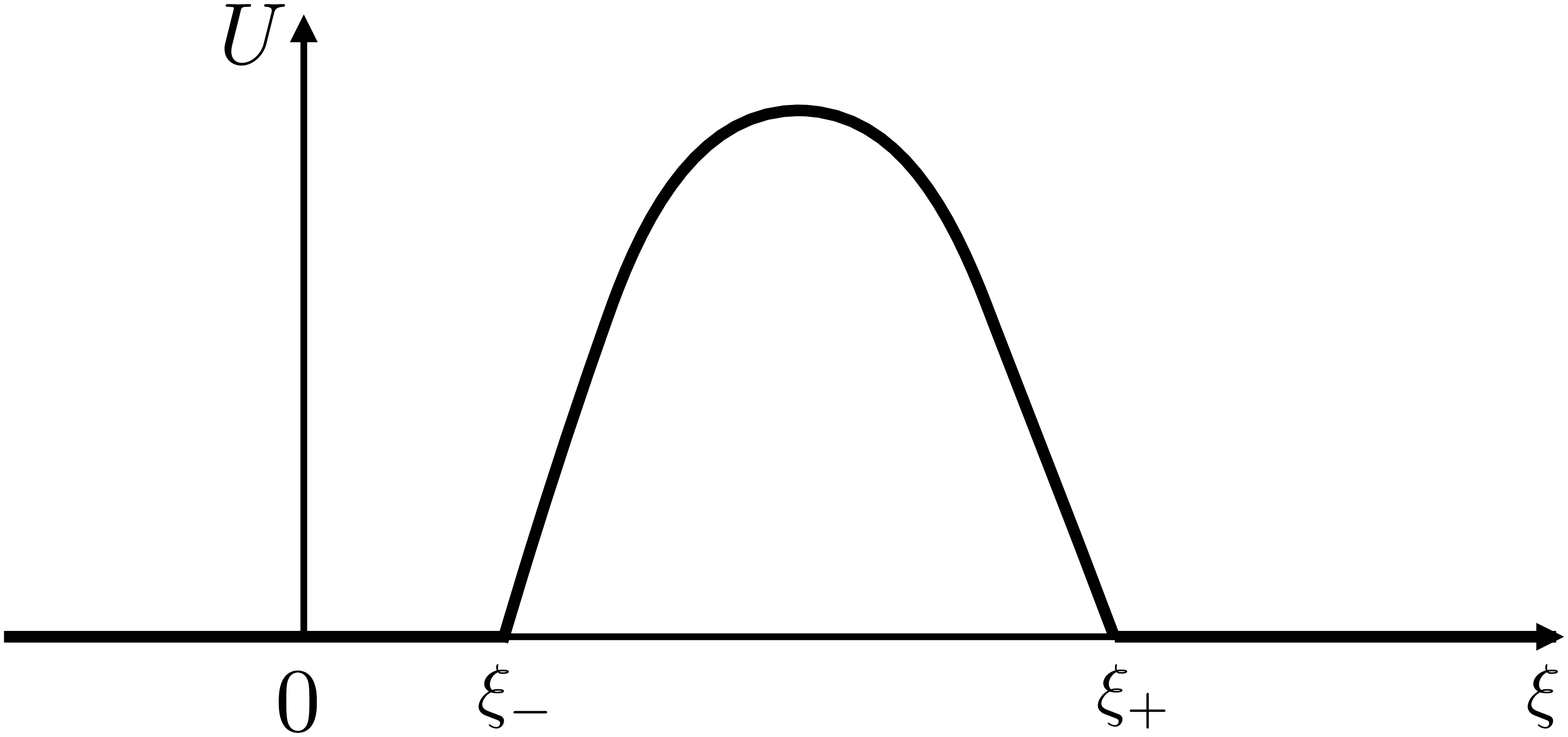}
\includegraphics[width=4.6cm]{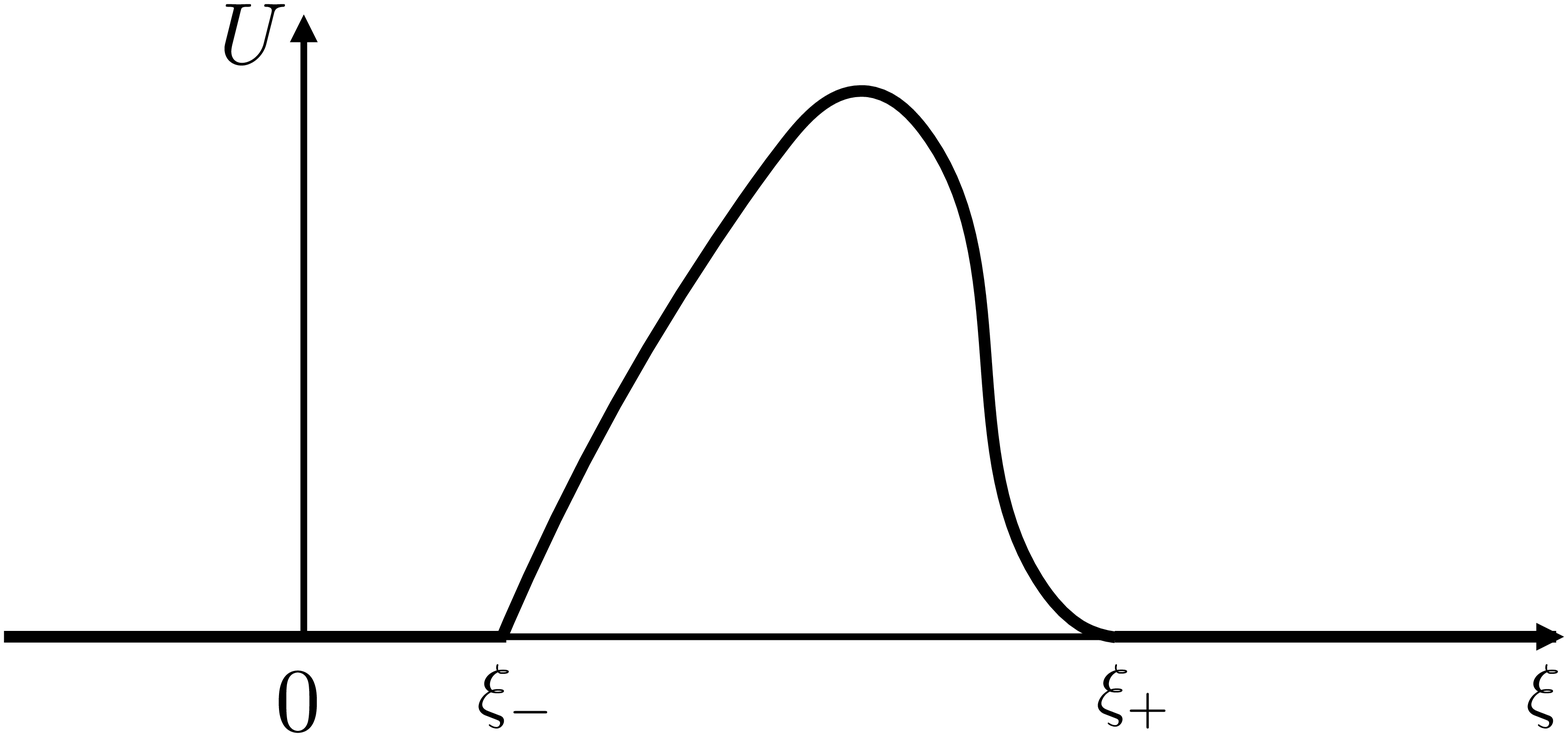}
\caption{Schematic pictures of the traveling wave solutions obtained in Theorem \ref{thm:pDNPE2-1}. Here it should be noted that the position of the singularity points $\xi_{-}$ and $\xi_{+}$ are not determined in our studies, however, they are shown in the figures for convenience.
[Left: The weak traveling wave solution in Theorem \ref{thm:pDNPE2-1} (i).]
[Center: The weak traveling wave solution in Theorem \ref{thm:pDNPE2-1} (ii).]
[Right: The weak traveling wave solution in Theorem \ref{thm:pDNPE2-1} (iii).]}
\label{fig:pDNPE2-1}
\end{figure}

Then, for the results of Proposition \ref{prop:pDNPE2-1} and Proposition \ref{prop:pDNPE2-2}, we use the transformation \eqref{eq:pDNPE1-2}.
Then, results for the classification of nonnegative traveling wave solutions (including weak sense) of \eqref{eq:pDNPE1-1} with $0<p<1$ for both $\delta=0$ and $\delta=1$ are immediately obtained.
As already mentioned, this argument is the same as \cite{cDNPE}.
Thus, the following results are compared to the nonnegative traveling wave solutions (including the weak sense) of \eqref{eq:pDNPE1-1} for $1<p \in \mathbb{R}$ in Corollary 1, Corollary 2, and Corollary 3 in \cite{cDNPE} in the case that $1<p \in \mathbb{R}$.
Figure \ref{fig:pDNPE2-1} is an image of the profile of the traveling wave solution (including the weak sense) obtained in Theorem \ref{thm:pDNPE2-1}.

\begin{thm}
\label{thm:pDNPE2-1}
Assume that $0<p<1$, $\mu>0$, and $\delta=0$.
Then, for a given positive constant $c$, the equation \eqref{eq:pDNPE1-1} has following three type of weak traveling wave solutions:

\begin{enumerate}
\item[(i)] 
There exists a family of weak traveling wave solutions (with the singularity).
Each solution $U(\xi)$ satisfies the following:
\begin{enumerate}
\item[(i1)]
$\displaystyle \lim_{\xi \searrow \xi_{-} + 0} U(\xi) = \lim_{\xi \to +\infty} U(\xi) = \lim_{\xi \to +\infty} U'(\xi) =  0$,
$\displaystyle\lim_{\xi \searrow \xi_{-} + 0} U'(\xi) = A_{2}$ with a positive constant $A_{2}$.
\item[(i2)]
$U(\xi)>0$ holds for $\xi\in (\xi_{-},+\infty)$ and $U(\xi)=0$ holds for $\xi\in (-\infty, \xi_{-}]$.
\item[(i3)]
There exists a constant $\xi_{0}\in (\xi_{-}, +\infty)$ such that the following holds:
$U'(\xi)>0$ for $\xi\in (\xi_{-}, \xi_{0})$, $U'(\xi_{0})=0$ and $U'(\xi)<0$ for $\xi\in (\xi_{0}, +\infty)$.
\end{enumerate}
In addition, the asymptotic behavior of $U(\xi)$ and $U'(\xi)$ for $\xi \searrow \xi_{-}+0$ are
\begin{equation}
\begin{cases}
U(\xi) \sim A_{1}(\xi-\xi_{-}),
\\
U'(\xi) \sim A_{2}
\end{cases}
\quad {\rm{as}} \quad \xi \searrow \xi_{-}+0,
\label{eq:pDNPE2-10}
\end{equation}
where $A_{j}$ are positive constants, and the asymptotic behavior of $U(\xi)$ and $U'(\xi)$ for $\xi \to +\infty$ are
\begin{equation}
\begin{cases}
U(\xi) \sim \left(\dfrac{c}{k}\xi\right)^{-1/p},
\\
U'(\xi) \sim -\dfrac{1}{p}\left(\dfrac{k}{c}\right)^{-\frac{1}{p}}\xi^{-\frac{p+1}{p}},
\end{cases}
\quad {\rm{as}} \quad \xi\to +\infty.
\label{eq:pDNPE2-11}
\end{equation}

\item[(ii)] 
There exists a family of weak traveling wave solutions (with the singularity).
Each solution $U(\xi)$ satisfies the following:
\begin{enumerate}
\item[(ii1)]
$\displaystyle \lim_{\xi \searrow \xi_{-} + 0} U(\xi) = \lim_{\xi \nearrow \xi_{+}-0} U(\xi) =0$,
$\displaystyle \lim_{\xi \searrow \xi_{-}+0}  U'(\xi) =  A_{2}$,
$\displaystyle\lim_{\xi \nearrow \xi_{+} - 0} U'(\xi) = -A_{4}$ with positive constants $C_{j}$.
\item[(ii2)]
$U(\xi)>0$ holds for $\xi\in (\xi_{-},\xi_{+})$ and $U(\xi)=0$ holds for $\xi\in (-\infty, \xi_{-}] \cup [\xi_{+}, +\infty)$.
\item[(ii3)]
There exists a constant $\xi_{0}\in (\xi_{-}, \xi_{+})$ such that the following holds:
$U'(\xi)>0$ for $\xi\in (\xi_{-}, \xi_{0})$, $U'(\xi_{0})=0$ and $U'(\xi)<0$ for $\xi\in (\xi_{0}, \xi_{+})$.
\end{enumerate}
In addition, the asymptotic behavior of $U(\xi)$ and $U'(\xi)$ for $\xi \searrow \xi_{-}+0$ are expressed as \eqref{eq:pDNPE2-10}, and the asymptotic behavior of $U(\xi)$ and $U'(\xi)$ for $\xi \nearrow \xi_{+}-0$ are
\begin{equation}
\begin{cases}
U(\xi) \sim A_{3}(\xi_{+}-\xi),
\\
U'(\xi) \sim -A_{4}
\end{cases}
\quad {\rm{as}} \quad \xi \nearrow \xi_{+}-0,
\label{eq:pDNPE2-12}
\end{equation}
where $A_{j}$ are positive constants.

\item[(iii)] 
There exists a weak traveling wave solution (with the singularity).
The solution $U(\xi)$ satisfies (ii2), (ii3), and the following:
\begin{enumerate}
\item[(iii1)]
$\displaystyle \lim_{\xi \searrow \xi_{-} + 0} U(\xi) = \lim_{\xi \nearrow \xi_{+}-0} U(\xi) =0$,
$\displaystyle \lim_{\xi \nearrow \xi_{+}-0}  U'(\xi) =  +\infty$,
$\displaystyle\lim_{\xi \searrow \xi_{-} + 0} U'(\xi) = -C$ with $C>0$. 
\end{enumerate}
In addition, the asymptotic behavior of $U(\xi)$ and $U'(\xi)$ for $\xi \searrow \xi_{-}+0$ are expressed as \eqref{eq:pDNPE2-10}, and the asymptotic behavior of $U(\xi)$ and $U'(\xi)$ for $\xi \nearrow \xi_{+}-0$ are
\begin{equation}
U(\xi) \sim A_{1}(\xi_{+}-\xi)^{\frac{1}{p}},
\quad
U'(\xi) \sim -A_{2}(\xi_{+}-\xi)^{-\frac{p-1}{p}}
\quad {\rm{as}} \quad \xi \nearrow \xi_{+}-0
\label{eq:pDNPE2-13}
\end{equation}
with positive constants $A_{j}>0$.
\end{enumerate}
\end{thm}

\begin{rem}
\label{rem:pDNPE2-1}
The result in (i) of Theorem \ref{thm:pDNPE2-1} corresponds to Corollary 1 of \cite{cDNPE}, only the asymptotic behavior in $\xi \searrow \xi_{-}+0$ is different.
On the other hand, \eqref{eq:pDNPE2-11} in Theorem \ref{thm:pDNPE2-1} agrees with the result of \cite{cDNPE}.

By comparing the results in \cite{cDNPE} with those in Theorem \ref{thm:pDNPE2-1}, results (ii) and (iii) suggest that traveling wave solutions that do not appear for $1<p \in \mathbb{R}$ do appear for $0<p<1$.
In other words, it suggests that the structure of traveling wave solutions changes at $p=1$.
In $p=1$, if the same Poincar\'e compactification method as in this paper is employed to investigate the structure of the traveling wave solution, there is a possibility that the center manifold exists near the equilibrium at infinity of the two-dimensional ODEs derived from the traveling wave coordinates.
The approximation of the center manifold and the flow on it will be studied, and the information will be used to determine the asymptotic behavior of the original traveling wave solution.
However, it is currently not possible to obtain a good approximation of the center manifold, as in \cite{DNLA, DNPE}, which is necessary to obtain the asymptotic behavior.
See Remark \ref{rem:pDNPE-U2-1}.
Therefore, the case $p=1$ requires more careful analysis and should be discussed separately from this paper.
\end{rem}

\begin{rem}
\label{rem:XuYin1}
Note that the profile of traveling waves around $\xi_{-}$ in Theorem \ref{thm:pDNPE2-1} corresponds to the non-$C^{1}$ sharp type in \cite{XuYin} (in detail, piecewise-$C^{1}$ sharp waves). 
\end{rem}

\begin{thm}
\label{thm:pDNPE2-2}
Assume that $0<p<1$, $\mu>0$, and $\delta=1$.
Then, for a given positive constant $c$, the equation \eqref{eq:pDNPE1-1} has three types of weak traveling wave solutions and one traveling wave solution including (ii) and (iii) in Theorem  \ref{thm:pDNPE2-1}.
The remaining two characterizations are as follows:

\begin{enumerate}
\item[(iv)] 
There exists a family of weak traveling wave solutions (with the singularity).
Each solution $U(\xi)$ satisfies (i2), (i3), and the following:
\begin{enumerate}
\item[(iv1)]
$\displaystyle 
\lim_{\xi \searrow \xi_{-} + 0} U(\xi)=0$,
$\displaystyle
\lim_{\xi \to +\infty} U(\xi) = \mu^{-1/p}$,
$\displaystyle
\lim_{\xi \searrow \xi_{-} + 0} U'(\xi) = C$ with a positive constant $C$.
\end{enumerate}
In addition, the asymptotic behavior of $U(\xi)$ and $U'(\xi)$ for $\xi \searrow \xi_{-}+0$ are expressed as \eqref{eq:pDNPE2-10}, and the asymptotic behavior of $U(\xi)$ for $\xi \to +\infty$ is
\begin{equation}
U(\xi) \sim \mu^{-\frac{1}{p}} \sim
\begin{cases}
\left(B_{1}e^{\omega_{1}\xi}+B_{2}e^{\omega_{2}\xi}+\dfrac{1}{\mu}\right)^{\frac{1}{p}}, \quad (D>0),
\\
\left((B_{3}\xi+B_{4})e^{\omega\xi}+\dfrac{1}{\mu}\right)^{\frac{1}{p}},\quad (D=0),
\\
\left(e^{-\frac{\mu c}{2}\xi} \bar{Z}(\xi) +\dfrac{1}{\mu}\right)^{\frac{1}{p}},\quad (D<0),
\end{cases}
\label{eq:pDNPE2-14}
\end{equation}
where $B_{j}$ ($1\le j\le 4$) are constants.

\item[(v)] 
There exists a traveling wave solution such that it corresponds to the orbit in \eqref{eq:pDNPE1-1}.
The solution $U(\xi)$ satisfies the following:
\begin{enumerate}
\item[(v1)]
$\displaystyle 
\lim_{\xi \to -\infty} U(\xi) =0,
\lim_{\xi \to +\infty} U(\xi) = \mu^{-1/p}.
$
\item[(v2)]
$U(\xi)>0$ holds for $\xi\in \mathbb{R}$.
\end{enumerate}
In addition, the asymptotic behavior of $U(\xi)$ for $\xi \to \infty$ is expressed as \eqref{eq:pDNPE2-14}, and the asymptotic behavior of $U(\xi)$ for $\xi \to -\infty$ is
\begin{equation}
U(\xi) \sim \left(\dfrac{Mc^{2}e^{\frac{p}{c}\xi}}{M(\mu c^{2}+1)e^{\frac{p}{c}\xi}-1}\right)^{\frac{1}{p}}
\quad {\rm{as}}\quad \xi\to -\infty,
\label{eq:pDNPE2-15} 
\end{equation}
where $M<0$ is a constant that depends on the initial state $\phi(0)=\phi_{0}$.
\end{enumerate}
\end{thm}

\begin{rem}
\label{rem:pDNPE2-2}
The result in Theorem \ref{thm:pDNPE2-2}(iv) corresponds to Corollary 2 in \cite{cDNPE}, and \eqref{eq:pDNPE2-14} agrees with the result of \cite{cDNPE}.
Note that the asymptotic behavior in $\xi \searrow \xi_{-}+0 $ is different.
The result in Theorem \ref{thm:pDNPE2-2}(v) corresponds to Corollary 3 in \cite{cDNPE}, and \eqref{eq:pDNPE2-15} matches the result in \cite{cDNPE}.
By comparing the results in \cite{cDNPE} with those in Theorem \ref{thm:pDNPE2-2}, it suggests that the structure of the traveling wave solution changes around $p=1$, as in Remark \ref{rem:pDNPE2-1}.
However, the discussion on traveling wave solutions in the case $p=1$ is open as in Remark \ref{rem:pDNPE2-1}.
\end{rem}

\begin{rem}
\label{rem:XuYin2}
Note that the traveling wave profile in Theorem \ref{thm:pDNPE2-2}(iv) corresponds to non-monotone wavefront piecewise-$C^{1}$ type sharp waves in \cite{XuYin}.
\end{rem}

\begin{figure}[t]
\centering
\includegraphics[width=6cm]{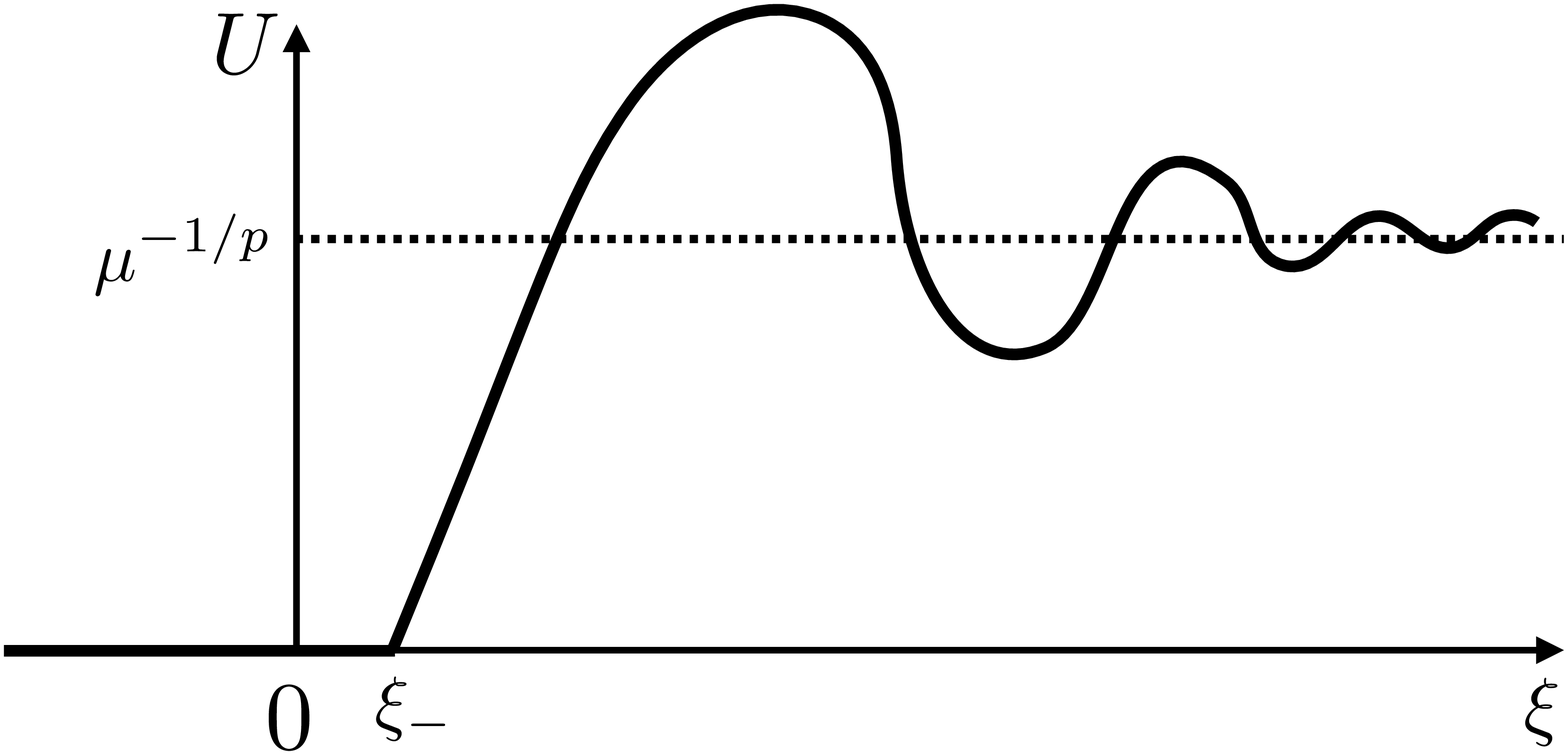}
\includegraphics[width=6cm]{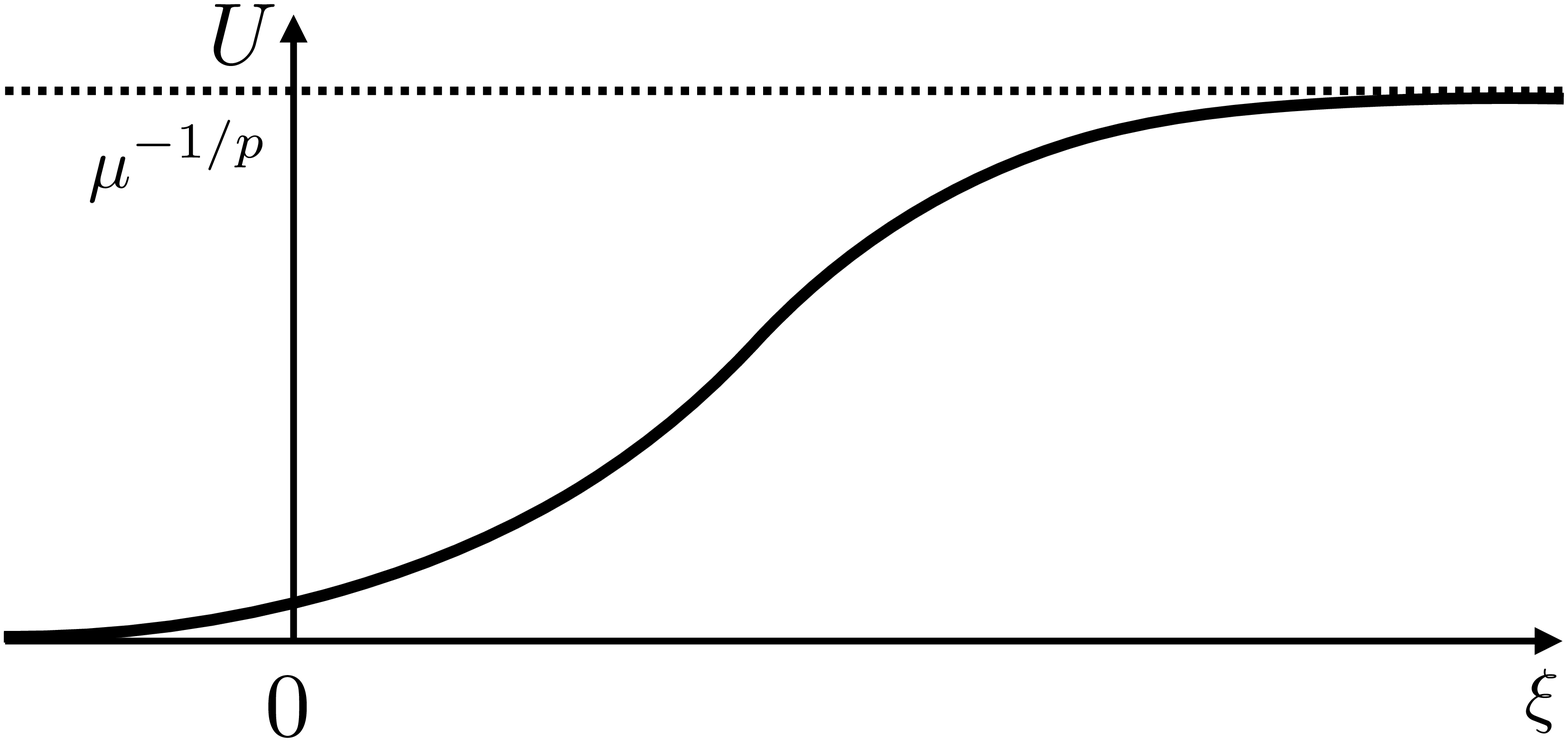}
\caption{Schematic pictures of the traveling wave solutions obtained in Theorem \ref{thm:pDNPE2-2}. Here it should be noted that the position of the singularity point $\xi_{-}$ is not determined in our studies, however, they are shown in the figures for convenience.
[Left: The weak traveling wave solution in Theorem \ref{thm:pDNPE2-2} (iv) in the case that $D<0$.]
[Right: The traveling wave solution on $\xi\in \mathbb{R}$ obtained in Theorem \ref{thm:pDNPE2-2} (v) in the case that $D>0$.]}
\label{fig:pDNPE2-2}
\end{figure}

Figure \ref{fig:pDNPE2-2} is an image of the profile of the traveling wave solution (including the weak sense) obtained in Theorem \ref{thm:pDNPE2-2}.

Next, the results of Theorem \ref{thm:pDNPE2-1} and Theorem \ref{thm:pDNPE2-2} and the transformation \eqref{eq:pDNPE1-4} are used.
Then, we immediately obtain results for the classification of nonnegative traveling wave solutions for $m>1$ (including weak sense) in \eqref{eq:pDNPE1-5} for both $\delta=0$ and $\delta=1$.

\begin{cor}
\label{cor:pDNPE2-1}
Assume that $m>1$, $k>0$, $\mu>0$, and $\delta=0$.
Then, for a given positive constant $c$, the equation \eqref{eq:pDNPE1-5} has following three type of weak traveling wave solutions:

\begin{enumerate}
\item[(A)] 
There exists a family of weak traveling wave solutions (with quenching).
Each solution $V(\xi)$ satisfies the following:
\begin{enumerate}
\item[(A1)]
$\displaystyle \lim_{\xi \searrow \xi_{-} + 0} V(\xi) = \lim_{\xi \to +\infty} V(\xi) = \lim_{\xi \to +\infty} V'(\xi) =  0$,
$\displaystyle\lim_{\xi \searrow \xi_{-} + 0} V'(\xi) = +\infty$.
\item[(A2)]
$V(\xi)>0$ holds for $\xi\in (\xi_{-},+\infty)$ and $V(\xi)=0$ holds for $\xi\in (-\infty, \xi_{-}]$.
\item[(A3)]
There exists a constant $\xi_{0}\in (\xi_{-}, +\infty)$ such that the following holds:
$V'(\xi)>0$ for $\xi\in (\xi_{-}, \xi_{0})$, $V'(\xi_{0})=0$ and $V'(\xi)<0$ for $\xi\in (\xi_{0}, +\infty)$.
\end{enumerate}
In addition, the asymptotic behavior of $V(\xi)$ and $V'(\xi)$ for $\xi \searrow \xi_{-}+0$ are
\begin{equation}
\begin{cases}
V(\xi) \sim A_{1}(\xi-\xi_{-})^{\frac{1}{m}},
\\
V'(\xi) \sim A_{2}(\xi-\xi_{-})^{-\frac{m-1}{m}}
\end{cases}
\quad {\rm{as}} \quad \xi \searrow \xi_{-}+0,
\label{eq:pDNPE2-16}
\end{equation}
where $A_{j}$ are positive constants, and the asymptotic behavior of $V(\xi)$ and $V'(\xi)$ for $\xi \to +\infty$ are
\begin{equation}
\begin{cases}
V(\xi) \sim m^{-\frac{1}{m-1}}\left(\dfrac{k}{c}\right)^{-\frac{1}{m-1}}\xi^{-\frac{1}{m-1}},
\\
V'(\xi) \sim -(m-1)^{-1}m^{-\frac{1}{m-1}}\left(\dfrac{k}{c}\right)^{-\frac{1}{m-1}}\xi^{-\frac{m}{m-1}},
\end{cases}
\quad {\rm{as}} \quad \xi\to +\infty.
\label{eq:pDNPE2-17}
\end{equation}

\item[(B)] 
There exists a family of weak traveling wave solutions (with quenching).
Each solution $V(\xi)$ satisfies the following:
\begin{enumerate}
\item[(B1)]
$\displaystyle \lim_{\xi \searrow \xi_{-} + 0} V(\xi) = \lim_{\xi \to \xi_{+}-0} V(\xi) =0$,
$\displaystyle \lim_{\xi \to \xi_{-}+0}  V'(\xi) =  +\infty$,
$\displaystyle\lim_{\xi \searrow \xi_{+} - 0} V'(\xi) = -\infty$.
\item[(B2)]
$V(\xi)>0$ holds for $\xi\in (\xi_{-},\xi_{+})$ and $V(\xi)=0$ holds for $\xi\in (-\infty, \xi_{-}] \cup [\xi_{+}, +\infty)$.
\item[(B3)]
There exists a constant $\xi_{0}\in (\xi_{-}, \xi_{+})$ such that the following holds:
$V'(\xi)>0$ for $\xi\in (\xi_{-}, \xi_{0})$, $V'(\xi_{0})=0$ and $V'(\xi)<0$ for $\xi\in (\xi_{0}, \xi_{+})$.
\end{enumerate}
In addition, the asymptotic behavior of $V(\xi)$ and $V'(\xi)$ for $\xi \searrow \xi_{-}+0$ are expressed as \eqref{eq:pDNPE2-16}, and the asymptotic behavior of $V(\xi)$ and $V'(\xi)$ for $\xi \nearrow \xi_{+}-0$ are
\begin{equation}
\begin{cases}
V(\xi) \sim A_{3}(\xi_{+}-\xi)^{\frac{1}{m}},
\\
V'(\xi) \sim -A_{4}(\xi_{+}-\xi)^{-\frac{m-1}{m}}
\end{cases}
\quad {\rm{as}} \quad \xi \nearrow \xi_{+}-0,
\label{eq:pDNPE2-18}
\end{equation}
where $A_{j}$ are positive constants.

\item[(C)] 
If $1<m<2$, then there exists a weak traveling wave solution (with quenching and singularity).
The solution $V(\xi)$ satisfies (B2), (B3), and the following:
\begin{enumerate}
\item[(C1)]
$\displaystyle \lim_{\xi \searrow \xi_{-} + 0} V(\xi) = \lim_{\xi \nearrow \xi_{+}-0} V(\xi) =0$,
$\displaystyle \lim_{\xi \searrow \xi_{-}+0}  V'(\xi) =  +\infty$,
$\displaystyle\lim_{\xi \nearrow \xi_{+} - 0} V'(\xi) = 0$. 
\end{enumerate}
In addition, the asymptotic behavior of $V(\xi)$ and $V'(\xi)$ for $\xi \searrow \xi_{-}+0$ are expressed as \eqref{eq:pDNPE2-16}, and the asymptotic behavior of $V(\xi)$ and $V'(\xi)$ for $\xi \nearrow \xi_{+}-0$ are
\begin{equation}
\begin{cases}
V(\xi) \sim A_{5}(\xi_{+}-\xi)^{\frac{1}{m-1}},
\\
V'(\xi) \sim -A_{6}(\xi_{+}-\xi)^{-\frac{m-2}{m-1}}
\end{cases}
\quad {\rm{as}} \quad \xi \nearrow \xi_{+}-0,
\label{eq:pDNPE2-19}
\end{equation}
where $A_{j}$ are positive constants.

If $m=2$, then there exists a weak traveling wave solution (with quenching and singularity).
The solution $V(\xi)$ satisfies (B2), (B3), and the following:
\begin{enumerate}
\item[(C2)]
$\displaystyle \lim_{\xi \searrow \xi_{-} + 0} V(\xi) = \lim_{\xi \nearrow \xi_{+}-0} V(\xi) =0$,
$\displaystyle \lim_{\xi \searrow \xi_{-}+0}  V'(\xi) =  +\infty$,
$\displaystyle\lim_{\xi \nearrow \xi_{+} - 0} V'(\xi) = -C$ with $C>0$. 
\end{enumerate}
In addition, the asymptotic behavior of $V(\xi)$ and $V'(\xi)$ for $\xi \searrow \xi_{-}+0$ are expressed as \eqref{eq:pDNPE2-16}, and the asymptotic behavior of $V(\xi)$ and $V'(\xi)$ for $\xi \nearrow \xi_{+}-0$ are expressed as \eqref{eq:pDNPE2-19}.

If $m>2$, then there exists a weak traveling wave solution (with quenching and singularity).
The solution $V(\xi)$ satisfies (B1), (B2) and (B3).
In addition, the asymptotic behavior of $V(\xi)$ and $V'(\xi)$ for $\xi \searrow \xi_{-}+0$ are expressed as \eqref{eq:pDNPE2-16}, and the asymptotic behavior of $V(\xi)$ and $V'(\xi)$ for $\xi \nearrow \xi_{+}-0$ are expressed as \eqref{eq:pDNPE2-19}.
\end{enumerate}
\end{cor}

\begin{cor}
\label{cor:pDNPE2-2}
Assume that $m>1$, $k>0$, $\mu>0$, and $\delta=0$.
Then, for a given positive constant $c$, the equation \eqref{eq:pDNPE1-5} has three types of weak traveling wave solutions and one traveling wave solution including (B) and (C) in Corollary \ref{cor:pDNPE2-1}.
The remaining two characterizations are as follows:

\begin{enumerate}
\item[(D)] 
There exists a family of weak traveling wave solutions (with quenching).
Each solution $V(\xi)$ satisfies (A2), (A3), and the following:
\begin{enumerate}
\item[(D1)]
$\displaystyle 
\lim_{\xi \searrow \xi_{-} + 0} V(\xi)=0$,
$\displaystyle
\lim_{\xi \to +\infty} V(\xi) = m^{-1/(m-1)}\mu^{-(m-1)/m^{2}}$,
$\displaystyle
\lim_{\xi \searrow \xi_{-} + 0} V'(\xi) = +\infty$.
\end{enumerate}
In addition, the asymptotic behavior of $V(\xi)$ and $V'(\xi)$ for $\xi \searrow \xi_{-}+0$ are expressed as \eqref{eq:pDNPE2-16}, and the asymptotic behavior of $V(\xi)$ for $\xi \to +\infty$ is
\begin{align}
V(\xi)
& \sim m^{-1/(m-1)}\mu^{-(m-1)/m^{2}}
\label{eq:pDNPE2-20}
\\
&\sim
\begin{cases}
m^{-\frac{1}{m-1}} \left( B_{1}e^{\omega_{1}\xi}+B_{2}e^{\omega_{2}\xi}+\dfrac{1}{\mu} \right)^{-(m-1)/m^{2}}, \quad (D>0),
\\
m^{-\frac{1}{m-1}} \left( (B_{3}\xi+B_{4})e^{\omega\xi}+\dfrac{1}{\mu} \right)^{-(m-1)/m^{2}},\quad (D=0),
\\
m^{-\frac{1}{m-1}} \left( e^{-\frac{\mu c}{2}\xi} \bar{Z}(\xi) +\dfrac{1}{\mu} \right)^{-(m-1)/m^{2}},\quad (D<0),
\end{cases}
\label{eq:pDNPE2-21}
\end{align}
where $B_{j}$ ($1\le j\le 4$) are constants.

\item[(E)] 
There exists a traveling wave solution.
The solution $V(\xi)$ satisfies the following:
\begin{enumerate}
\item[(E1)]
$\displaystyle 
\lim_{\xi \to -\infty} V(\xi) =0,
\lim_{\xi \to +\infty} V(\xi) = m^{-1/(m-1)}\mu^{-(m-1)/m^{2}}.
$
\item[(E2)]
$V(\xi)>0$ holds for $\xi\in \mathbb{R}$.
\end{enumerate}
In addition, the asymptotic behavior of $V(\xi)$ for $\xi \to \infty$ is expressed as \eqref{eq:pDNPE2-21}, and the asymptotic behavior of $V(\xi)$ for $\xi \to -\infty$ is
\begin{equation}
V(\xi) \sim 
m^{-\frac{1}{m-1}} \left( \dfrac{Mc^{2}e^{\frac{p}{c}\xi}}{M(\mu c^{2}+1)e^{\frac{p}{c}\xi}-1} \right)^{(m-1)/m^{2}}
\quad {\rm{as}}\quad \xi\to -\infty,
\label{eq:pDNPE2-22} 
\end{equation}
where $M<0$ is a constant that depends on the initial state $\phi(0)=\phi_{0}$.
\end{enumerate}
\end{cor}

Finally, we use the classification results in Corollary \ref{cor:pDNPE2-1} with $\delta=0$ and $\mu=0$. 
Thus, we obtain the classification results for nonnegative weak traveling wave solutions in the porous medium equation \eqref{eq:pDNPE1-po}.

\begin{cor}
\label{cor:pDNPE2-3}
Assume that $m>1$, $\mu=0$, and $\delta=0$.
Then, for a given positive constant $c$, the equation \eqref{eq:pDNPE1-po} has three types of weak traveling wave solutions same as Corollary \ref{cor:pDNPE2-1}.
(B) and (C) in Corollary \ref{cor:pDNPE2-1} also hold here.
However, since the asymptotic behavior \eqref{eq:pDNPE2-17} of the weak traveling wave solution in (A) (see Corollary \ref{cor:pDNPE2-1}) as $\xi\to +\infty$ vanishes at $\mu=0$ ($k=p\mu$, see \eqref{eq:pDNPE1-2}), we only know that (A1) holds.
\end{cor}

\begin{figure}
\centering
\includegraphics[width=6cm]{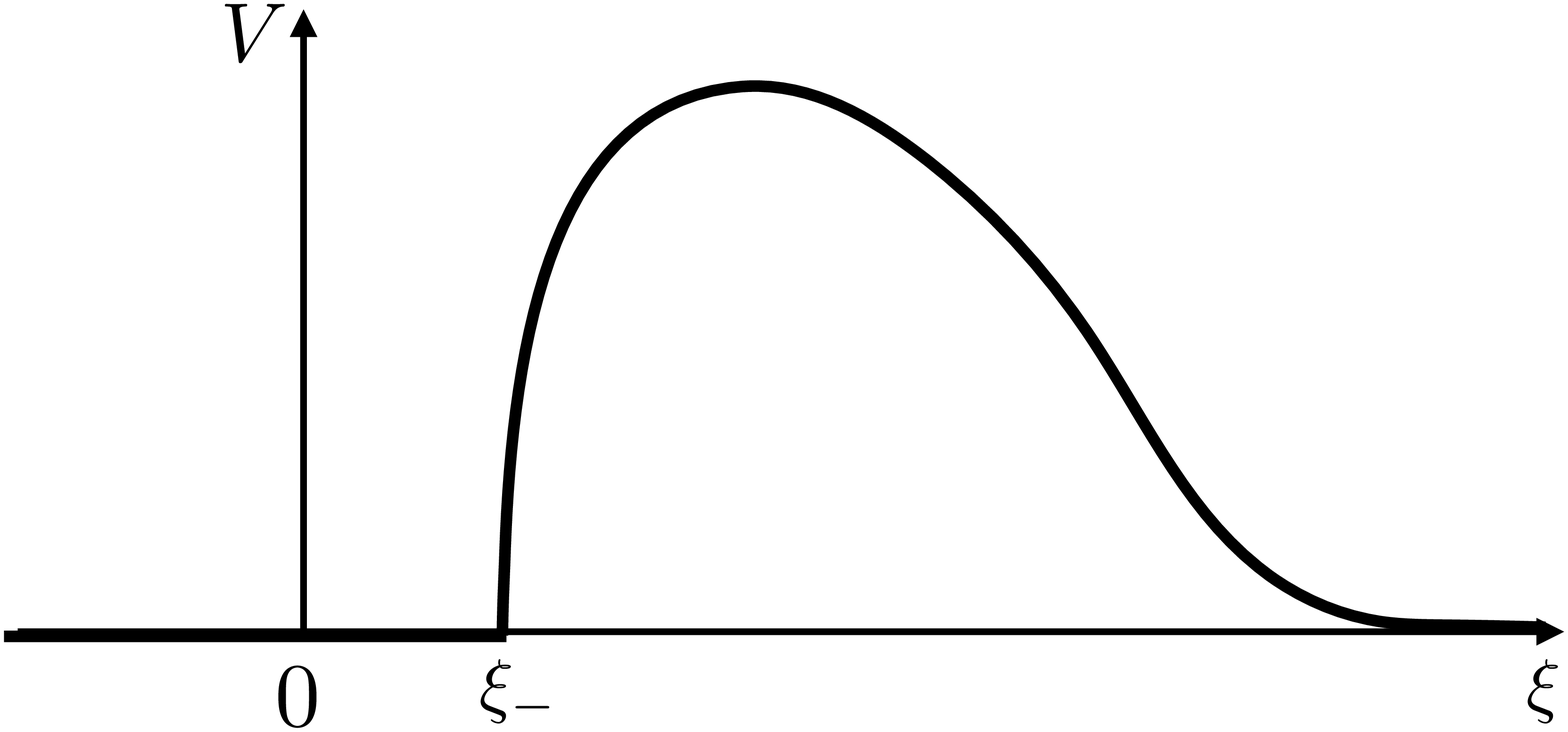}
\includegraphics[width=6cm]{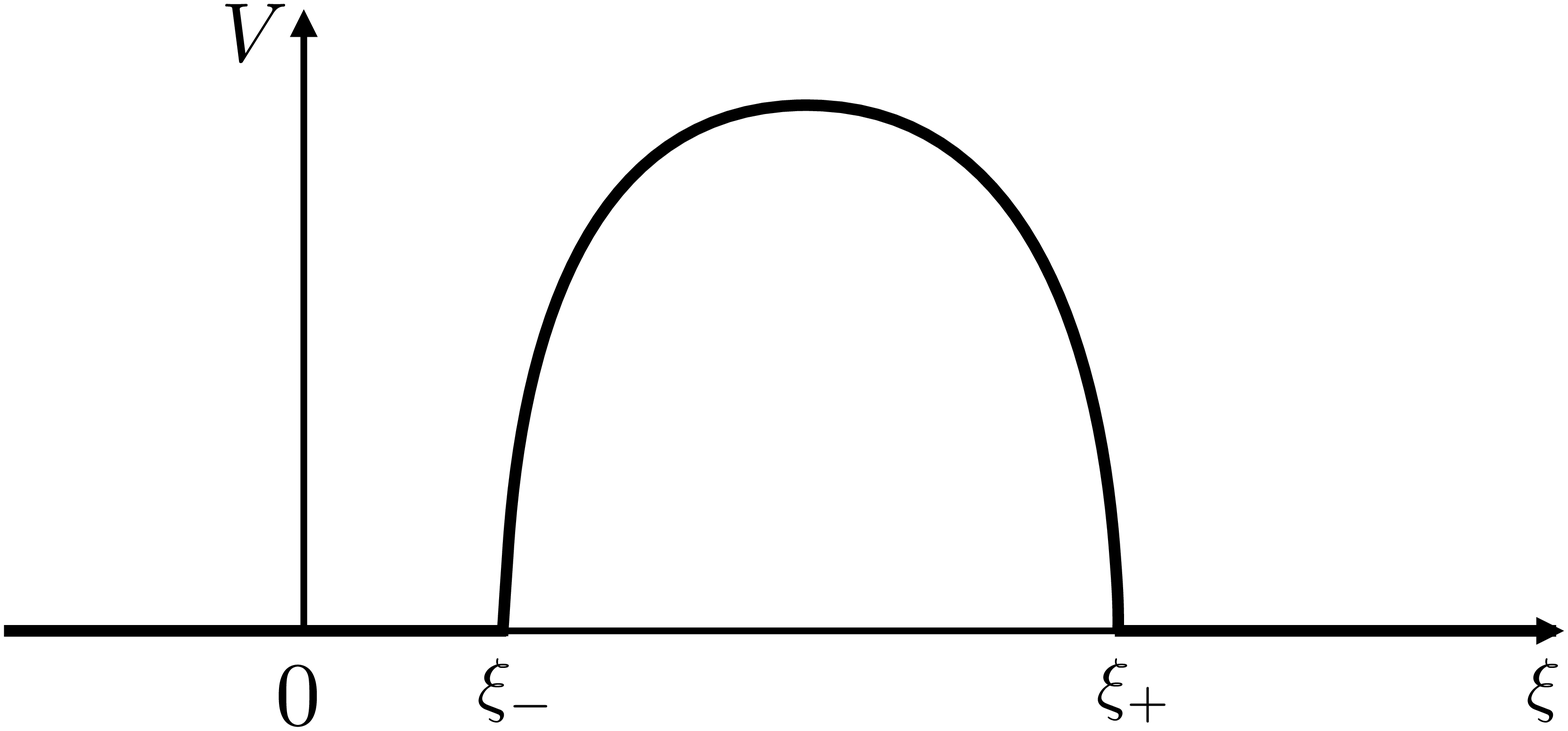}
\caption{Schematic pictures of the traveling wave solutions obtained in Corollary \ref{cor:pDNPE2-3}. Here it should be noted that the position of the singularity points $\xi_{-}$ and $\xi_{+}$ are not determined in our studies, however, they are shown in the figures for convenience.
[Left: The weak traveling wave solution in Corollary \ref{cor:pDNPE2-3} (A).]
[Right: The weak traveling wave solution in Corollary \ref{cor:pDNPE2-3} (B).]\\ \\
}
\label{fig:pDNPE2-3}
\end{figure}

\begin{figure}
\centering
\includegraphics[width=4.5cm]{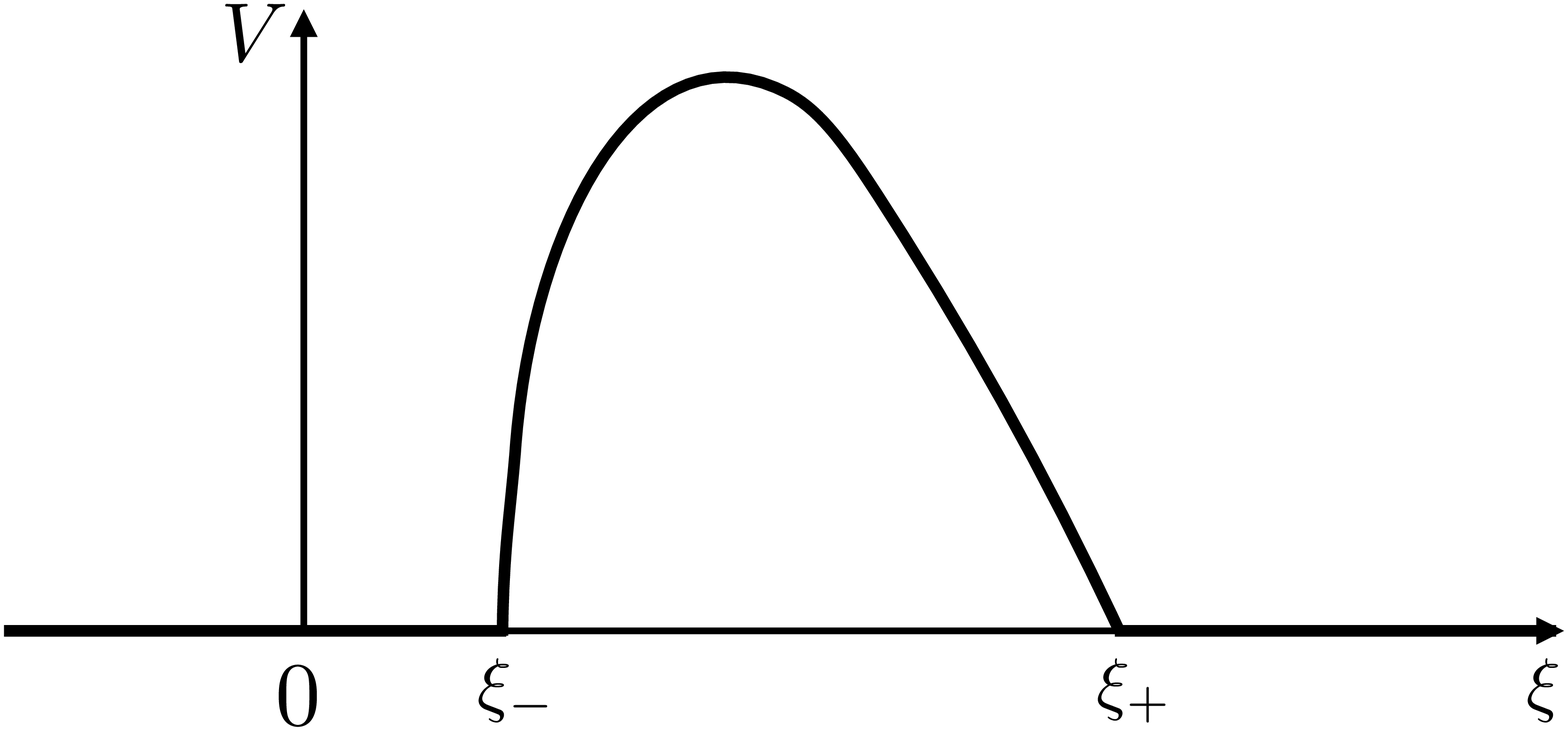}
\includegraphics[width=4.5cm]{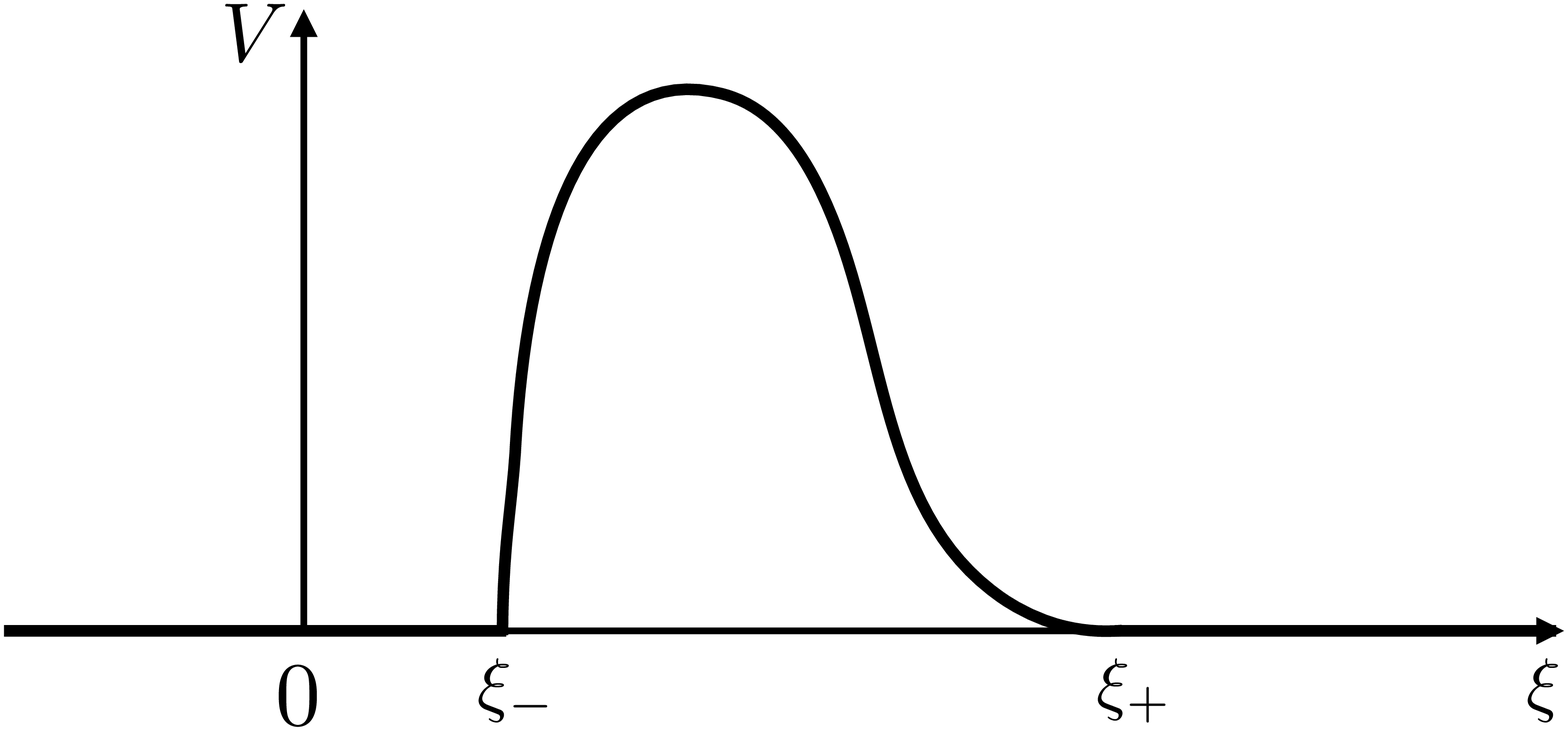}
\includegraphics[width=4.5cm]{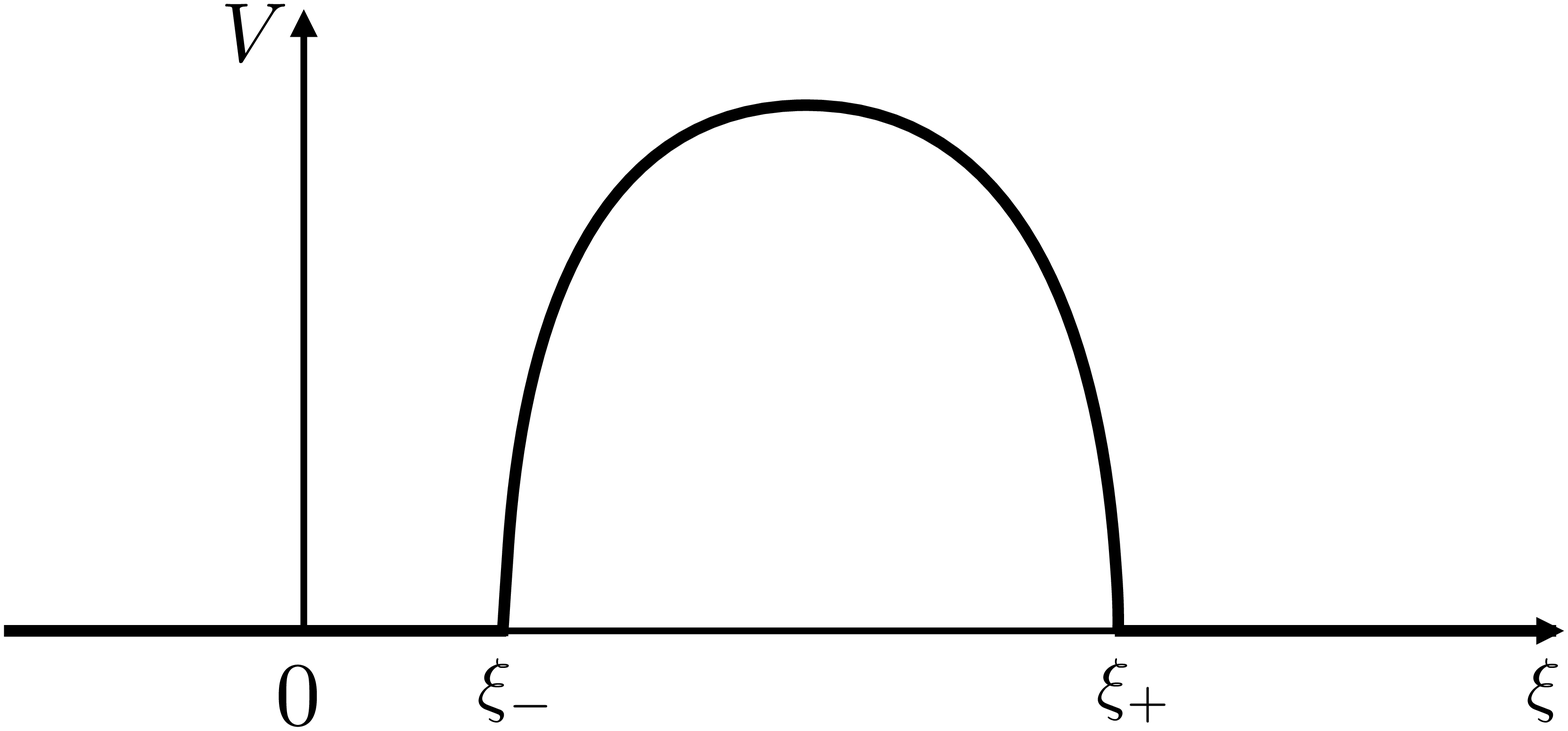}
\caption{Schematic pictures of the traveling wave solutions obtained in Corollary \ref{cor:pDNPE2-3}. Here it should be noted that the position of the singularity points $\xi_{-}$ and $\xi_{+}$ are not determined in our studies, however, they are shown in the figures for convenience.
[Left: The weak traveling wave solution in Corollary \ref{cor:pDNPE2-3} (C) in the case that $1<m<2$.]
[Center: The weak traveling wave solution in Corollary \ref{cor:pDNPE2-3} (C) in the case that $m=2$.]
[Right: The weak traveling wave solution in Corollary \ref{cor:pDNPE2-3} (C) in the case that $m>2$.]
}
\label{fig:pDNPE2-ad}
\end{figure}

Figure \ref{fig:pDNPE2-3} and Figure \ref{fig:pDNPE2-ad} is schematic pictures of the profile of the weak traveling wave solution obtained in Corollary \ref{cor:pDNPE2-3}.

\section{Dynamics on the Poincar\'e disk of \eqref{eq:pDNPE1-7}}
\label{sec:pDNPE3}
In this section, we study $\mathbb{R}^{2} \cup \{(\phi, \psi) \mid \|(\phi, \psi)\|=+\infty\}$ called the dynamics on the Poincar\'e disk, by the Poincar\'e compactification.
The discussion is almost the same as for \cite{cDNPE}, but for the reader's convenience, we briefly discuss its details and differences from \cite{cDNPE}.

\subsection{Dynamics near finite equilibria}
\label{sub:pDNPE3-1}
First, we study the dynamics near finite equilibria of \eqref{eq:pDNPE1-7}.
When $\delta=0$, there is no equilibrium.
If $\delta=1$, this equation has an equilibrium $E_{1} : (\phi,\psi)=(\mu^{-1},0)$ at $\{\phi>0\}$.
Note that $\phi=0$ has a singularity.
The existence of this equilibrium and the behavior around it is exactly the same as the discussion in \cite{cDNPE}.
This paper also uses the information around this equilibrium, so the results are described in this paper, although they overlap with those of \cite{cDNPE}.
The Jacobian matrix of the vector field \eqref{eq:pDNPE1-7} in $E_{1}$ is 
\[
E_{1}: \left(\begin{array}{cc}
0 & 1
\\
-k &-c\mu
\end{array}\right).
\]
Let $J_{1}$ be this matrix.
Then, the behavior of the solution around $E_{1}$ is different by the sign of $D=\mu^{2}c^{2}-4k$.
For instance, the matrix $J_{1}$ has the real distinct eigenvalues if $D>0$ and other cases can be concluded similarly.
If $c>0$, then $E_{1}$ is a stable node for $D\ge 0$, and is a stable focus (spiral sink) for $D<0$.

Second, in order to study the dynamics of \eqref{eq:pDNPE1-7} on the Poincar\'e disk, we desingularize it by the time-scale desingularization: 
\begin{equation} 
ds/d\xi = \phi^{-1}
\label{eq:pDNPE3-1}
\end{equation}
as in \cite{cDNPE}.
See Remark 9 in \cite{cDNPE} for a note on the time scale desingularization.
In addition, we refer to Section 7.7 of \cite{CK} and references therein for the analytical treatments of 
desingularization with time rescaling. 
Note that this allows us to include $\phi=0$.
Since we are considering a nonnegative solution, i.e., $\phi\ge 0$, the direction of the time does not change via this desingularization in this region.
Then we have
\begin{equation}
\begin{cases}
\phi'=\phi\psi,
\\
\psi'=-c\psi+\gamma\psi^{2}-k\phi^{2}+\delta p \phi,
\end{cases}
\quad \left( \, '=\dfrac{d}{ds}\right).
\label{eq:pDNPE3-2}
\end{equation}

\eqref{eq:pDNPE3-2} has the following equilibria regardless of $\delta$:
\[
E_{0}: (\phi, \psi)=(0,0), \quad E_{2}: (\phi, \psi)=\left( 0, c\gamma^{-1}\right).
\]
Note that it differs from \cite{cDNPE} in that $\gamma<0$.
The Jacobian matrices of the vector field \eqref{eq:pDNPE3-2} at these equilibria are
\[
E_{0}: \left(\begin{array}{cc}
0 & 0\\
\delta p &-c
\end{array}\right),
\quad
E_{2}: \left(\begin{array}{cc}
c\gamma^{-1} & 0
\\
\delta p &c
\end{array}\right).
\]
\begin{itemize}
\item
When $\delta=0$, $E_{2}$ is a saddle since the eigenvalues are $c \gamma^{-1}<0$ and $c>0$.
The eigenvectors for each eigenvalue are $(1,0)^{T}$ and $(0,1)^{T}$ with $T$ representing the transpose.
Then, the center manifold theory is applicable to study the dynamics near $E_{0}$ as in \cite{cDNPE} (for instance, see \cite{carr}). 
We can obtain the approximation of the (graph of) center manifold as follows:
\begin{equation}
\{ (\phi, \psi) \mid \psi(s)=-kc^{-1}\phi^{2}+O(\phi^{4})  \}. 
\label{eq:pDNPE3-3} 
\end{equation}
Hence, the dynamics of \eqref{eq:pDNPE3-2} near $E_{0}$ is topologically equivalent to the dynamics of the following equation:
\begin{equation}
\phi'(s)=-kc^{-1}\phi^{3}+O(\phi^{5}).
\label{eq:pDNPE3-4} 
\end{equation}
These results are the same as in Subsection 3.1 of \cite{cDNPE}.
See Subsection 3.2 of \cite{DNPE} for a more detailed process.
\item
When $\delta=1$, $E_{2}$ is a saddle.
The eigenvectors for each eigenvalue are $(1, c^{-1}p(\gamma^{-1}-1)^{-1})^{T}$ and $(0,1)^{T}$
As noted in \cite{cDNPE}, by the same argument as in Subsection 3.2 of \cite{DNPE} and Section 2 of \cite{DNLA}, we can study the dynamics around the equilibrium $E_{0}$ for $\delta=1$.
We conclude that the approximation of the (graph of) center manifold is
\begin{equation}
\left\{ (\phi, \psi) \mid  \psi(s)=pc^{-1}\phi -pc^{-3}(\mu c^{2}+1)\phi^{2}+O(\phi^{3}) \right\}
\label{eq:pDNPE3-5}
\end{equation}
and the dynamics of \eqref{eq:pDNPE3-2} near $E_{0}$ is topologically equivalent to the dynamics of the following equation:
\begin{equation}
\phi'(s)=pc^{-1}\phi^{2}-pc^{-3}(\mu c^{2}+1)\phi^{3}+O(\phi^{4}).
\label{eq:pDNPE3-6}
\end{equation}
\end{itemize}

Since the signs of $\gamma$ are only different in \cite{cDNPE}, the location of the finite equilibrium $E_{2}$ changes and stability and dynamics near $E_{0}$ and $E_{2}$ remain unchanged.
From the next subsection, we can consider the dynamics of this equation on the charts $\overline{U}_{j}$ ($j=1,2$) and $\overline{V}_{2}$ since we are considering $\phi\ge 0$.

\subsection{Dynamics on the chart $\overline{U}_{2}$}
\label{sub:pDNPE3-2}
To obtain the dynamics on the chart $\overline{U}_{2}$, we introduce the coordinates $(\lambda, x)$ by the formulas
\[
\phi(s)=x(s)/\lambda(s), \quad \psi(s)=1/\lambda(s).
\]
In this chart, it corresponds to $\phi\to 0$ and $\psi\to +\infty$ and the direction in which $x$ is positive corresponds to the direction in which $\phi$ is positive.
For a geometric image, see Fig. 2 of \cite{QTW} and Fig. 2 of \cite{DNPE}.
Then, we have
\begin{equation}
\begin{cases}
\lambda'=c\lambda -\gamma+kx^{2}-\delta p\lambda x,
\\
x'=p^{-1}\lambda^{-1}x+cx+k\lambda^{-1}x^{3}-\delta px^{2}.
\end{cases}
\label{eq:pDNPE3-7}
\end{equation}
By using the time-scale desingularization $d\tau/ds=\lambda^{-1}$, we can obtain
\begin{equation}
\begin{cases}
\lambda_{\tau}=c\lambda^{2}-\gamma \lambda+k\lambda x^{2}-\delta p\lambda^{2}x,
\\
x_{\tau}=p^{-1}x+c\lambda x+kx^{3}-\delta p\lambda x^{2},
\end{cases}
\label{eq:pDNPE3-8}
\end{equation}
where $\lambda_{\tau}=d\lambda/d\tau$ and $x_{\tau}=dx/d\tau$.
Up to this point, it is the same as \cite{cDNPE}, but the discussion of the following equilibrium is different.
The equilibrium of the system \eqref{eq:pDNPE3-8} on $\{\lambda=0\}$ is
\[
E_{3}: (\lambda,x)=(0,0).
\]
The Jacobian matrix of the vector field \eqref{eq:pDNPE3-8} at this equilibrium is
\[
E_{3}: \left(\begin{array}{cc}
-\gamma & 0\\
0 &1/p
\end{array}\right).
\]
Therefore, $E_{3}$ is source.
Note that $p>1$ in \cite{cDNPE} is a saddle since $\gamma$ is $0<\gamma<1$.

\begin{rem}
\label{rem:pDNPE-U2-1}
When $p=1$, the dynamics near $E_{3}$ can be investigated by applying the center manifold theorem.
This point requires careful and hard analysis as described in Remark \ref{rem:pDNPE2-1} and is an open problem.
\end{rem}

\subsection{Dynamics on the charts $\overline{V}_{2}$ and $\overline{U}_{1}$}
\label{sub:pDNPE3-3}
By introducing the same transformations and time scale transformations as in cite{cDNPE}, we obtain the same equations in each local coordinate.
Similarly to Subsection \ref{sub:pDNPE3-2}, we conclude that the equilibrium $E_{4}: (\lambda,x)=(0,0)$ exists and is a sink.
On the other hand, the differential equation system in $\overline{U}_{1}$ has no equilibrium, as in the argument about \cite{cDNPE}.

\subsection{Dynamics and connecting orbits on the Poincar\'e disk}
\label{sub:pDNPE3-5}
\begin{figure}[t]
\centering
\includegraphics[width=5cm]{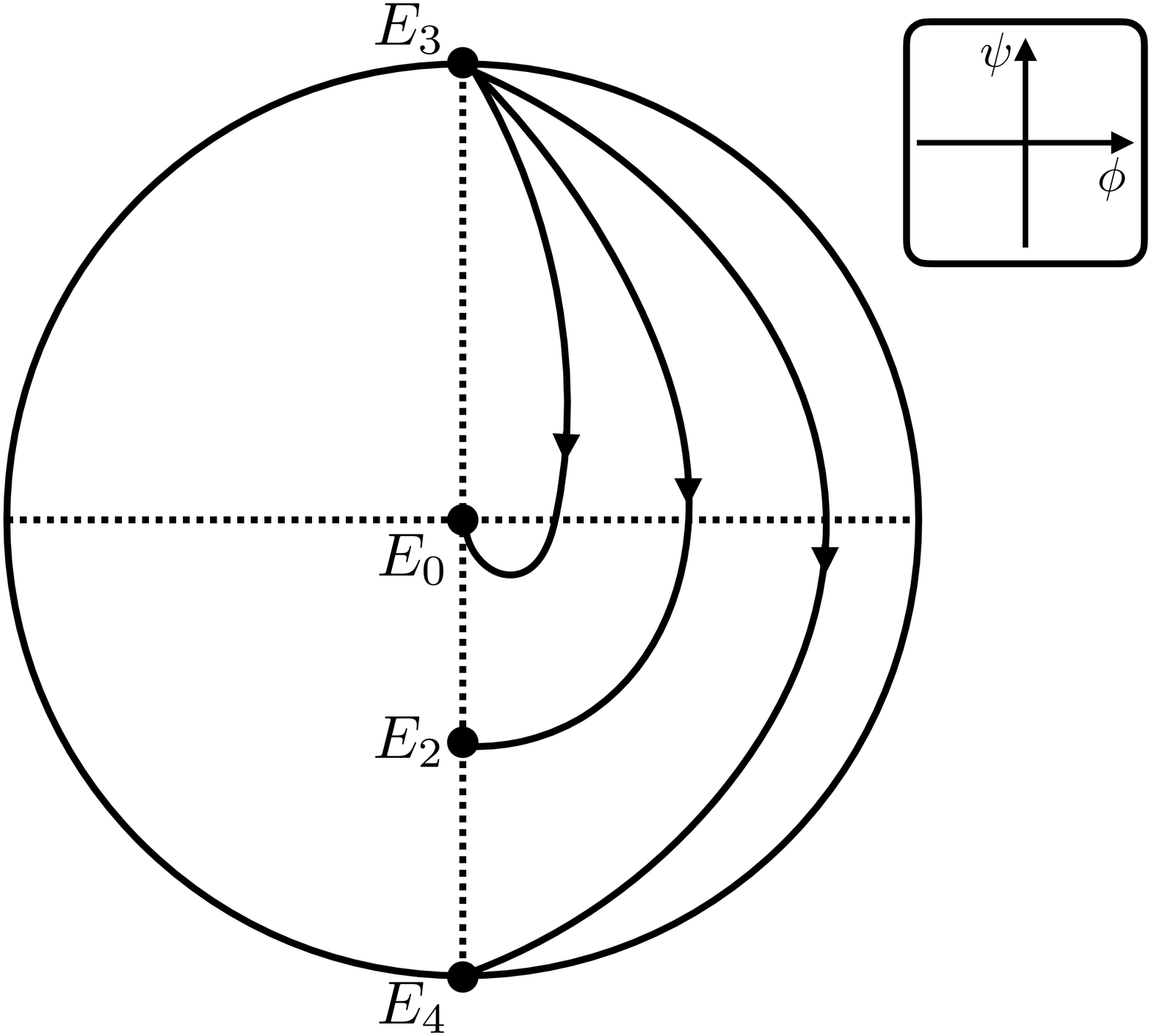}
\includegraphics[width=5cm]{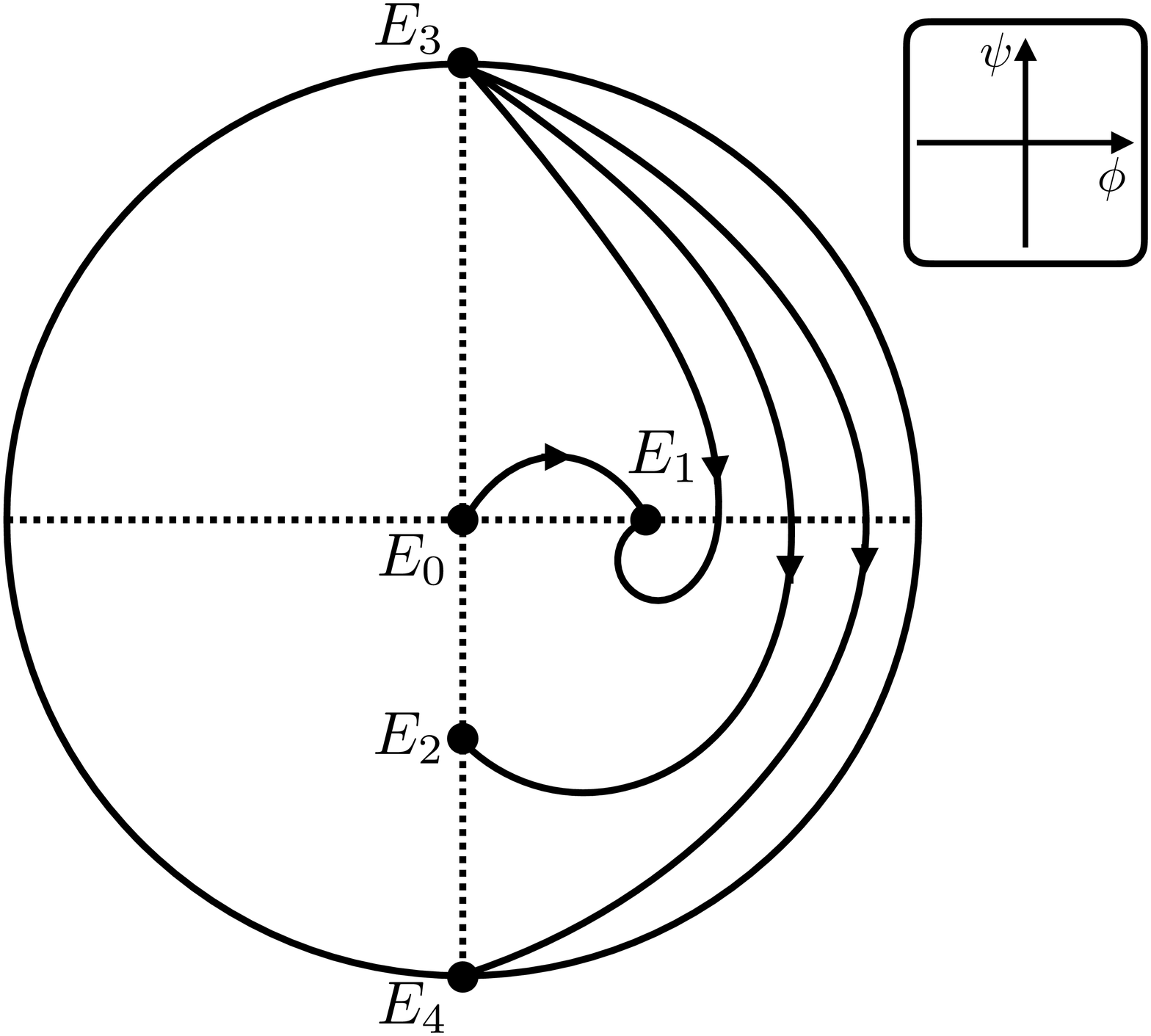}
\caption{Schematic pictures of the dynamics on the Poincar\'e disk in the case that  $0<p<1$.  
[Left: $\delta=0$.]
[Right: $\delta=1$.]
}
\label{fig:pDNPE3-1}
\end{figure}

Combining the dynamics on the charts $\overline{U}_{j}$ ($j=1,2$) and $\overline{V}_{2}$, we can obtain the dynamics on the Poincar\'e disk that is equivalent to the dynamics of \eqref{eq:pDNPE1-7} (or \eqref{eq:pDNPE3-2}) (see also Figure \ref{fig:pDNPE3-1}).
The set $\Phi$ denotes $\Phi=\{ (\phi, \psi) \in \mathbb{R}^{2} \cup \{ \|(\phi, \psi)\|=+\infty \} \}$.

The purpose of this subsection is to prove the existence of connecting orbits (see also Figure \ref{fig:pDNPE3-1}). 
The proof strategy for connecting orbits is the same as for \cite{DNPE, DNLA, cDNPE}.
Note, however, that the position of the equilibria is more complicated than in \cite{DNPE, DNLA, cDNPE}.
Before we do so, we will give some remarks about disks.

\begin{rem}[\cite{cDNPE}, Remark 11]
\label{rem:pDNPE3-2}
In Figure \ref{fig:pDNPE3-1}, we need to be careful about the handling of the point $E_{0}$.
A note on this treatment is given for the reader's convenience, although it is a reproduction of \cite{DNPE, DNLA, cDNPE}.
When we consider the parameter $s$ on the disk, $E_{0}$ is the equilibrium of \eqref{eq:pDNPE3-2}.
However, $E_{0}$ is a point on the line $\{\phi=0\}$ with singularity about the parameter $\xi$.
We see that $d \phi/d \psi$ takes the same values on the vector fields defined by \eqref{eq:pDNPE1-7} and \eqref{eq:pDNPE3-2} except the singularity $\{\phi=0\}$.
If the trajectories start (resp. come in) the equilibrium $E_{0}$ about the parameter $s$, then they start from (resp. come in) the point $E_{0}$ about $\xi$. 
This is also the case for $E_{2}$, $E_{3}$, and $E_{4}$.
\end{rem}

First, we show the existence of connecting orbits at $\delta=0$ (see also Figure \ref{fig:pDNPE3-1}).
Since $\{\phi=0\}$ is invariant in \eqref{eq:pDNPE3-2}, any trajectory starting from a point on $\{(\phi, \psi)\in \Phi \mid \phi>0\}$ must not go to $\{(\phi, \psi) \in \Phi \mid \phi<0\}$.
Let $\mathcal{W}^{s}(E_{2})$ be a stable manifold of $E_{2}$ in the dynamical system \eqref{eq:pDNPE3-2} and $\mathcal{W}^{u}(E_{3})$ be an unstable manifold of $E_{3}$ in the dynamical system \eqref{eq:pDNPE3-2}
Then, the only trajectories that reach a point on $\mathcal{W}^{s}(E_{2})$ are those starting from a point on $\mathcal{W}^{s}(E_{3})$ from the Poincar\'e-Bendixson theorem.
That is, in \eqref{eq:pDNPE3-2}, the existence of a trajectory between $E_{3}$ and $E_{2}$ such that starting from a point on $\mathcal{W}^{u}(E_{3})$ and reaching a point on $\mathcal{W}^{s}(E_{2})$ is shown.
Here, as mentioned in Remark \ref{rem:pDNPE3-2}, since $d \phi/d \psi$ takes the same value on the vector fields defined by \eqref{eq:pDNPE1-7} and \eqref{eq:pDNPE3-2}, these trajectories obtained by \eqref{eq:pDNPE3-2} are also inherited by \eqref{eq:pDNPE1-7}.
Thus, in \eqref{eq:pDNPE1-7}, the existence of a trajectory between $E_{3}$ and $E_{2}$ such that starting from a point on $\mathcal{W}^{u}(E_{3})$ and reaching a point on $\mathcal{W}^{s}(E_{2})$ is shown.
Similarly, the existence of connecting orbits between $E_{3}$ and $E_{0}$ and between $E_{3}$ and $E_{4}$ at $\delta=0$ is also shown.

Next, we show the existence of connecting orbits for $\delta=1$.
$\mathcal{W}^{cu}(E_{0})$ denotes the center-unstable manifold of $E_{0}$ in the dynamical system \eqref{eq:pDNPE3-2} for $\delta=1$.
$\mathcal{W}^{s}(E_{1})$ denotes the stable manifold of $E_{1}$ in the dynamical system \eqref{eq:pDNPE1-7} and \eqref{eq:pDNPE3-2} for $\delta=1$.
$\mathcal{W}^{s}(E_{4})$ denotes the stable manifold of $E_{4}$ in the dynamical system \eqref{eq:pDNPE3-2} for $\delta=1$.
By applying the Poincar\'e-Bendixson theorem, a trajectory starting from a point on $\mathcal{W}^{cu}(E_{0})$ can only go to a point on $\mathcal{W}^{s}(E_{1})$ or $\mathcal{W}^{s}(E_{4})$.
Let us assume that the trajectory starting from a point on $\mathcal{W}^{cu}(E_{0})$ goes to $\mathcal{W}^{s}(E_{4})$.
Since
\[
\left. \psi' \vphantom{\big|}\right|_{0<\phi<\mu^{-1},\, \psi=0} =-p\phi(\mu\phi-1)>0
\]
in \eqref{eq:pDNPE3-2}, this trajectory will pass through a point on the $\phi$ axis satisfying $\phi>\mu^{-1}$ and head for a point on $\mathcal{W}^{s}(E_{4})$.
In this case, the trajectory toward a point on $\mathcal{W}^{s}(E_{1})$ is a trajectory starting from a point on $\mathcal{W}^{u}(E_{3})$.
However, this trajectory intersects the assumed trajectory, which is inconsistent.
The same argument leads to a contradiction if we assume that the trajectory starting from a point on $\mathcal{W}^{cu}(E_{0})$ heads onto $\mathcal{W}^{s}(E_{2})$.
From the above, we conclude that a trajectory starting from a point on $\mathcal{W}^{cu}(E_{0})$ can only go to a point on $\mathcal{W}^{s}(E_{1})$.
That is, in \eqref{eq:pDNPE3-2}, the existence of a connecting orbit between $E_{0}$ and $E_{1}$ is shown.
The same argument as for $\delta=0$ shows the existence of a connecting orbit between $E_{0}$ and $E_{1}$, $E_{3}$ and $E_{1}$, and $E_{3}$ and $E_{4}$ in the case that $\delta$ is $\delta=1$.

From the above, as shown in Figure \ref{fig:pDNPE3-1}, the existence of all connecting orbits on the Poincar\'e disk for both $\delta=0$ and $\delta=1$ is shown.

\section{Proof of the main results}
\label{sec:pDNPE4}
In this chapter, we prove the main results.
As mentioned in Section \ref{sec:pDNPE2}, we first prove Proposition \ref{prop:pDNPE2-1} and Proposition \ref{prop:pDNPE2-2}.
Then, using these results and \eqref{eq:pDNPE1-2}, from the same argument as \cite{cDNPE}, we obtain Theorem \ref{thm:pDNPE2-1} and Corollary \ref{cor:pDNPE2-1}, Theorem \ref{thm:pDNPE2-2} and Corollary \ref{cor:pDNPE2-2} are obtained.
Thus, the proofs of Proposition \ref{prop:pDNPE2-1} and Proposition \ref{prop:pDNPE2-2} are central and essential in this chapter.

If the initial data are located on $\Phi \backslash \{\phi< 0\}$, the existence of the solutions follows from the standard theory for the ODEs.
Therefore, we consider the existence of the trajectories that connect equilibria and detailed dynamics near the equilibria on the Poincar\'e disk and their asymptotic behavior.
Note from the argument of Subsection \ref{sub:pDNPE3-5} that there are no orbits that pass through the $\psi$-axis, so there are no sign-changing solutions.

\subsection{Proof of Proposition \ref{prop:pDNPE2-1}}
\label{sub:pDNPE4-1}
The discussion of Subsection \ref{sub:pDNPE3-5} reveals all connecting orbits of \eqref{eq:pDNPE1-7} with $\phi \ge 0$ for $\delta=0$.
Thus, we obtain the existence of three types of traveling waves corresponding to these three types: trajectories between $E_{3}$ and $E_{0}$, between $E_{3}$ and $E_{2}$, and between $E_{3}$ and $E_{4}$.
Note that there is only one connecting orbit between $E_{3}$ and $E_{2}$, but the other connecting orbits form a family of connecting orbits from the discussion of dynamical systems near the equilibria.

First, we show (I).
By the same argument as for \cite{DNPE, cDNPE}, we obtain (I2), (I3) and \eqref{eq:pDNPE2-5}, except that
\[
\lim_{\xi\to +\infty}u(\xi) =\lim_{\xi\to +\infty}u'(\xi)=0
\]
and $\xi_{-}$ exists.
For (I2), it follows from the fact that $\phi(\xi)>0$ holds that $u(\xi)>0$ holds.
For (I3), it is shown that the orbit passes through the $\phi$-axis, corresponding to $u'(\xi)=0$.
What else should be shown is the existence of $-\infty <\xi_{-}<+\infty$ and \eqref{eq:pDNPE2-4}.

The following discussion seems similar to that of \cite{cDNPE}. 
However, the difference is that \cite{cDNPE} uses an approximation of solutions near a finite equilibrium, while the following discussion uses an approximation of solutions near equilibrium at infinity.
It is written in the same way as \cite{cDNPE}, but for the reader's convenience, the proof is given below.

The solutions around $E_{3}$ are approximated as 
\[
\begin{cases}
\lambda(\tau)= A_{1}e^{-\gamma \tau}(1+o(1)),
\\
x(\tau)= A_{2}e^{p^{-1}\tau}(1+o(1))
\end{cases}
\]
with constants $A_{j}>0$ ($j=1,2$).
Using this equation, we obtain
\begin{align*}
\dfrac{d\tau}{d\xi}
&= \dfrac{d\tau}{ds}\dfrac{ds}{d\xi}
=\lambda^{-1}\phi^{-1}
=x^{-1}
\sim A_{3}e^{-p^{-1}\tau} 
\quad {\rm{as}}\quad \tau\to -\infty.
\end{align*}
From this result, we can obtain $d\xi/d\tau \sim A_{4}e^{p^{-1}\tau}$.
This yields $\xi(\tau)\sim A_{5}e^{p^{-1}\tau}+A_{6}$.
However, $A_{6}$ is a constant and the other $A_{j}$ is a positive constant.
Set 
\[
\xi_{-}=\lim_{\tau\to -\infty}\xi(\tau),
\]
then, we have 
\[
\xi-\xi_{-}\sim A_{7}e^{p^{-1}\tau} 
\quad {\rm{as}}\quad \tau\to -\infty
\]
with a positive constant $A_{7}>0$ in the same discussion as \cite{QTW, DNPE, cDNPE, Matsue1, Matsue2}.
Therefore, we obtain 
\begin{align*}
\phi(\xi)
&= \lambda^{-1}x
\sim A_{8} e^{\tau}
\sim A_{9} (\xi-\xi_{-})^{p}
\quad {\rm{as}} \quad \xi \searrow \xi_{-}+0,
\\
\psi(\xi)
&=\lambda^{-1} 
\sim A_{10} e^{\gamma \tau}
\sim A_{11}(\xi-\xi_{-})^{p-1}
\quad {\rm{as}} \quad \xi \searrow \xi_{-}+0
\end{align*}
with positive constants $A_{j}$.
Since the trajectories are lying on $\{(\phi, \psi) \mid \phi>0, \psi>0\}$, it holds that $A_{9}>0$ and $A_{11}>0$ (see also Figure \ref{fig:pDNPE3-1}).
Since $u(\xi)=\phi(\xi)$ and $u_{\xi}=\psi(\xi)$ hold, we can derive \eqref{eq:pDNPE2-4}.

The same arguments as in (I) and \cite{DNPE, cDNPE} above can be used to obtain the results in (II) and (III).
By combining the above results, Proposition \ref{prop:pDNPE2-1} is proved.
\qed
\\

\subsection{Proof of Proposition \ref{prop:pDNPE2-2}}
\label{sub:pDNPE4-2}
Conclusions for $\delta=1$ can be obtained as in Subsection \ref{sub:pDNPE4-1} and \cite{DNPE, DNLA, cDNPE}.
See these references for details.
\qed
\\

\section{Discussion}
\label{sec:pDNPE5}
In this paper, we classify all connecting orbits on the phase space $\mathbb{R}^{2} \cup \{ \|(\phi, \psi)\|=+\infty \}, \backslash\, \{\phi<0\}$ of \eqref{eq:pDNPE1-7} using Poincar\'e compactification and classical dynamical systems theory. 
These characterize the non-negative (including in a weak sense) traveling wave solutions of \eqref{eq:pDNPE1-3}, whose existence, profile and asymptotic behavior are given by Proposition \ref{prop:pDNPE2-1} and Proposition \ref{prop:pDNPE2-2}.
Using these results and \eqref{eq:pDNPE1-2}, the classification of non-negative (including in a weak sense) traveling wave solutions for $0<p<1$ in \eqref{eq:pDNPE1-1} was obtained with Theorem \ref{thm:pDNPE2-1} and Theorem \ref{thm:pDNPE2-2}.
Note that this idea is essentially the same as \cite{cDNPE}, but with results in the $p$ range not obtained in the authors' previous works \cite{DNLA, DNPE}.
Furthermore, in this paper, the transformation of \eqref{eq:pDNPE1-4} allows us to obtain classification results for nonnegative weak traveling wave solutions of spatial $1$-dimensional porous medium equation such as Corollary \ref{cor:pDNPE2-3}.

In this section, we treat all the dynamics of \eqref{eq:pDNPE1-7} obtained in \cite{cDNPE} and Section \ref{sec:pDNPE3}, including at infinity.
We describe the stability changes of the equilibria at infinity and finite equilibria, i.e., the bifurcation of the equilibria in \eqref{eq:pDNPE1-7} with p as a parameter.

\begin{figure}[t]
\centering
\includegraphics[width=14cm]{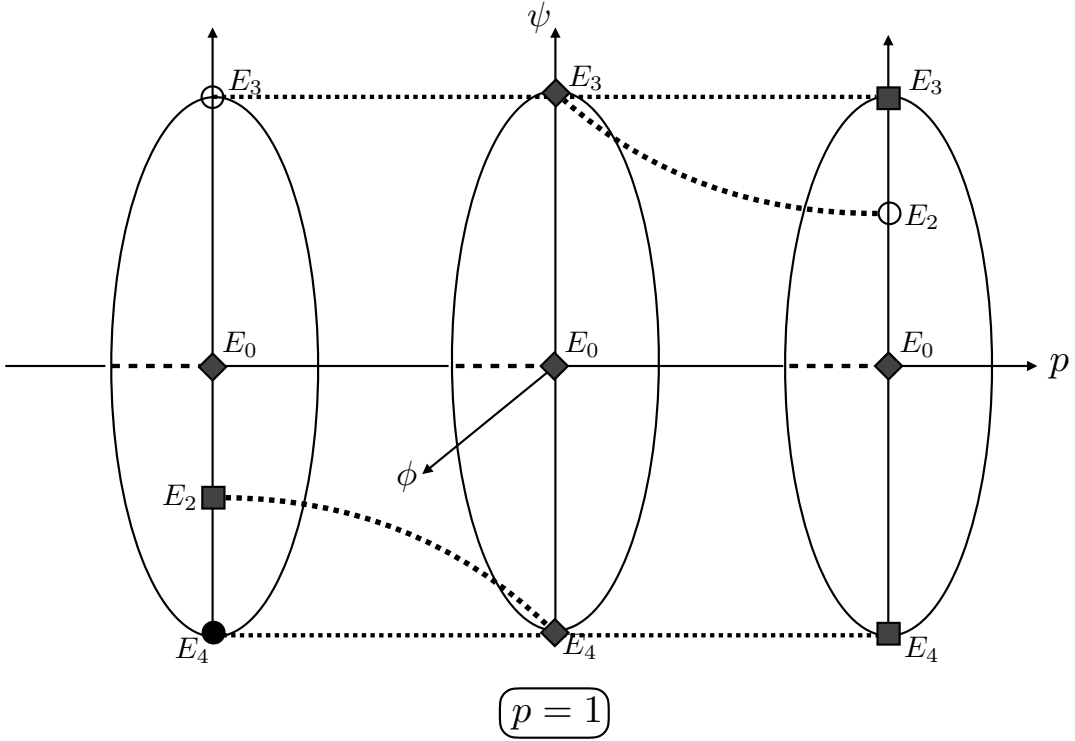}
\caption{
Schematic diagram of the finite equilibria and equilibria at infinity on the  Poincar\'e disk. 
White and Black circles denote source and sink equilibria, respectively.  Gray squares denote saddle equilibria, and gray diamonds equilibria with the center manifold.
Left, center and right disks correspond to the case that  $0<p<1$, $p=1$ and $p>1$, respectively.  
}
\label{fig:pDNPE5-1}
\end{figure}

In Figure \ref{fig:pDNPE5-1}, the positions of equilibria on the $\psi$-axis on the Poincar\'e disk for $0<p<1$ and $p=1$ and $p>1$ and their respective stability are shown using symbols.
Since $E_{3}$ and $E_{4}$ are equilibrium points at infinity, the $\psi$-axis in Figure \ref{fig:pDNPE5-1} represents
\[
\{-\infty\} \cup \{ -\infty <\psi< +\infty\} \cup \{+\infty\}.
\]
See the caption of Figure \ref{fig:pDNPE5-1} for the notation of the stability of the equilibria on each $p$.
By comparing the results for $p>1$ obtained in \cite{cDNPE} with the results for $0<p<1$ obtained in this paper, we can observe the change in stability with respect to finite equilibria and equilibria at infinity in \eqref{eq:pDNPE1-7}.
Moving $p$ from $p>1$ to $0<p<1$ with $p=1$ as a critical value, the change in stability of each equilibrium can be summarized as follows: 
\begin{itemize}
\item 
The equilibrium at infinity $E_{3}: (\phi, \psi)=(0, +\infty)$ is a saddle at $p>1$, but becomes a source at $0<p<1$ after $p=1$.
That is, the stable and unstable dimensions, which were $1$ each for $p>1$, lose $1$ of stable dimension at $p=1$, and for $0<p<1$, the stable dimension is $0$ and the unstable dimension is $2$.
\item 
The equilibrium point at infinity $E_{4}: (\phi, \psi)=(0, -\infty)$ is similar to $E_{3}$.
\item 
For a finite equilibrium $E_{0}:(\phi, \psi)=(0,0)$ is an equilibrium for which the center manifold theorem applies with $0<p<1$ and $p=1$ and $p>1$.
\item
For the finite equilibrium $E_{2}: (\phi, \psi)=(0, pc/(p-1))$, the unstable dimension is $2$ when $p>1$.
It exists in $\psi<0$ when $0<p<1$ and the unstable dimension is $1$.
When $p=1$, $E_{2}$ does not exist.
This can be interpreted as $\psi=pc/(p-1)\to \pm\infty$ with $p \to 1\pm 0$ (corresponding to the gray dotted line in Figure \ref{fig:pDNPE5-1}).
\end{itemize}

\begin{figure}[t]
\centering
\includegraphics[width=11cm]{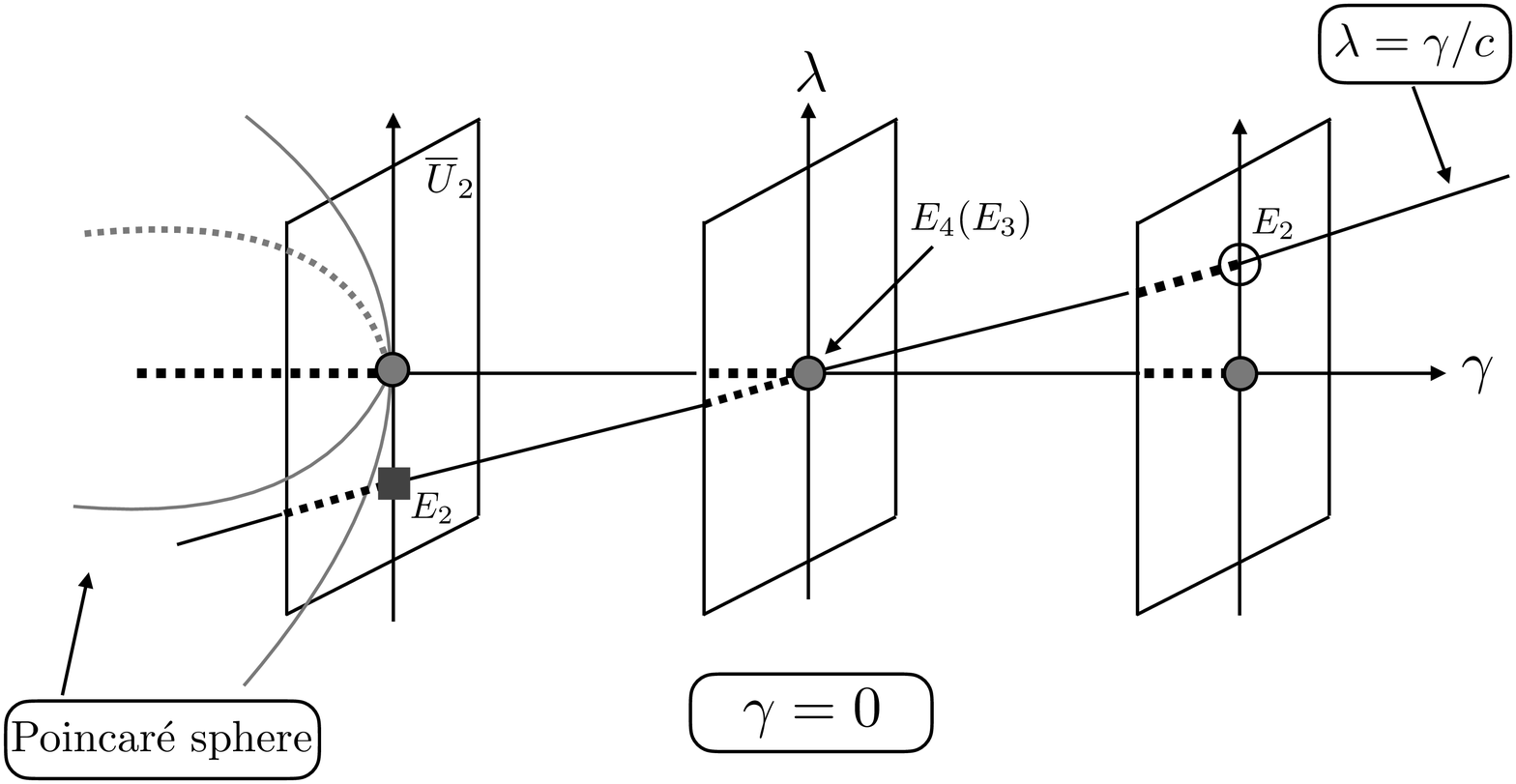}
\caption{
Schematic picture of the bifurcation on $\overline{U}_2$. 
The gray circles correspond to the equilibria $E_3$ and $E_4$ on $\overline{U}_2$. 
}
\label{fig:bif1}
\end{figure} 

Moreover, it follows that for a given $\varepsilon>0$, there exist the dynamical systems that are not topologically equivalent on the region $\{ |p -1| < \varepsilon\}$. 
Indeed, the dynamics of \eqref{eq:pDNPE3-8} restricted on $\{x=0\}$ is 
\begin{equation}
\lambda_{\tau} = (c \lambda -\gamma) \lambda, \quad (\gamma = (p-1)/p, \, c>0). \label{eq:bif}
\end{equation}
Then $\gamma = (p-1)/p = (p-1) +O((p-1)^2)$ holds near $p=1$, therefore 
we can see that the trans-critical bifurcation occur at $\gamma=(p-1)/p=0$ in \eqref{eq:bif}, that is, the equilibria $\lambda=0$ and $\lambda=\gamma/c$ of \eqref{eq:bif} through the  trans-critical bifurcation. 

Since we consider the dynamics on $\{\phi>0\}$ of \eqref{eq:pDNPE1-7}, here we ignore the direction of the trajectories on $\psi$-axis and focus only the graph of bifurcated branches corresponding to the equilibria $\lambda=0$ and $\lambda = \gamma/c$ of \eqref{eq:bif} (see Figure \ref{fig:bif1}). 
 Then the dynamics on the local chart $\overline{U}_2$ and $\overline{V}_2$ are located as upside-down, and the bifurcated branches on the upper half plane of $\overline{U}_2$ and lower half plane of that are a map from $\{(\psi,\phi) \,|\, \psi >0\}$ and $\{(\psi,\phi) \,|\, \psi <0\}$ by the mapping $g_{2}^{+} \circ f^{+}$ and $g_2^{+} \circ f^{-}$, respectively (see Appendix A). 
Therefore, the bifurcated branches 
\[ 
\{(x,\lambda,\gamma) \,| \, (x,\lambda)=(0,0)\} \quad \mbox{and} \quad 
\{(x,\lambda,\gamma)\,|\, (x,\lambda) = (0,\gamma /c)\} 
\] 
on $\overline{U}_2$ are divided into two parts and observed on the Poincar\'e disk.
These branches correspond to the bifurcated branches 
\[ 
\{(\phi,\psi,p) \,|\, (\phi,\psi) = (0,pc/(p-1)), \, 0<p<1\} \cup \{(\phi,\psi)=(0,-\infty)\}
\]
and
\[ 
\{(\phi,\psi,p) \,|\, (\phi,\psi) = (0,pc/(p-1)), \, 1<p\} \cup \{(\phi,\psi)=(0,+\infty)\} 
\]
on the phase space of \eqref{eq:pDNPE1-7}, i.e., the bifurcations at infinity occur at $p=1$ in the dynamical system \eqref{eq:pDNPE1-7}.  
This implies that the traveling waves of \eqref{eq:pDNPE1-1} are characterized by the bifurcations at infinity observed in \eqref{eq:pDNPE1-7}.

This result is based on the investigation of all dynamics, including $\{ \|(\phi,\psi)\| = \infty \}$, of the two-dimensional ODEs \eqref{eq:pDNPE1-7} that characterize the traveling waves of \eqref{eq:pDNPE1-1} by Poincar\'e compactification.
To the best of the authors' knowledge, our results will provide a new perspective on the bifurcation analysis as well the as applications of Poincar\'e compactification to the partial differential equations.

\section*{Acknowledgments}
IY was partially supported by JSPS KAKENHI Grant Number JP21J20035.


\appendix
\section*{Appendix A: Overview of the Poincar\'e type compactification}
\label{appendix:ap1}
\setcounter{figure}{0}
\renewcommand{\thefigure}{A.\arabic{figure}}
The Poincar\'e compactification is one of the compactifications of the original phase space (the embedding of $\mathbb{R}^{n}$ into the unit upper hemisphere of $\mathbb{R}^{n+1}$).
In this appendix, we briefly introduce the Poincar\'e compactification. 
Here Section 2 of \cite{QTW, DNPE} are reproduced.
Also, it should be noted that we refer \cite{FAL} for more details.
Let 
\[ X = P(\phi,\psi) \dfrac{\partial \mbox{}}{\partial \phi} 
+ Q(\phi,\psi) \dfrac{\partial \mbox{}}{\partial \psi}\]
be a polynomial vector field on $\mathbb{R}^{2}$, or in other words
\[
\begin{cases}
\dot{\phi} = P(\phi,\psi),
\\
\dot{\psi} = Q(\phi,\psi),
\end{cases}
\]
where $\dot{~\mbox{}~}$ denotes $d/dt$, and $P$, $Q$ are polynomials of arbitrary degree in the variables $\phi$ and $\psi$.

We consider $\mathbb{R}^{2}$ as the plane in $\mathbb{R}^{3}$ defined by $(y_{1},y_{2},y_{3})=(\phi,\psi,1)$.
We consider the sphere $\mathbb{S}^{2} = \{ y \in \mathbb{R}^{3} \, |\, y_{1}^{2} + y_{2}^{2}+y_{3}^{2}=1\}$ which we call Poincar\'e sphere.
We divide the sphere into
\[
H_{+} = \{ y \in \mathbb{S}^{2}\,|\,y_{3}>0\},
\quad
H_{-} = \{ y \in \mathbb{S}^{2}\,|\,y_{3}<0\} 
\]
and
\[\mathbb{S}^{1} = \{y \in \mathbb{S}^{2}\, | \, y_{3}=0\}.\]

Let us consider the embedding of vector field $X$ from $\mathbb{R}^{2}$ to $\mathbb{S}^{2}$ given by
\[ f^{+}:\mathbb{R}^{2} \to \mathbb{S}^{2},
 \quad f^{-}:\mathbb{R}^{2} \to \mathbb{S}^{2},\]
where
\[f^{\pm}(\phi,\psi):= \pm \left( \dfrac{\phi}{\Delta(\phi,\psi)},\dfrac{\psi}{\Delta(\phi,\psi)},\dfrac{1}{\Delta(\phi,\psi)} \right) \]
with $\Delta(\phi,\psi) = \sqrt{\phi^{2}+\psi^{2}+1}$.

Then we consider six local charts on $\mathbb{S}^{2}$ given by $U_{k} = \{y \in \mathbb{S}^{2} \, | \, y_{k}>0\}$, $V_{k} = \{y \in \mathbb{S}^{2} \, | \, y_{k}<0\}$ for $k=1,2,3$.
Consider the local projection
\[ g^{+}_{k} : U_{k} \to \mathbb{R}^{2}, 
\quad g^{-}_{k} : V_{k} \to \mathbb{R}^{2} \]
defined as
\[g^{+}_{k}(y_{1},y_{2},y_{3}) = - g^{-}_{k}(y_{1},y_{2},y_{3})
 = \left(\dfrac{y_{m}}{y_{k}},\dfrac{y_{n}}{y_{k}} \right) \]
 for $m<n$ and $m,n \not= k$. 
 The projected vector fields are obtained as the vector fields on the planes
\[
\overline{U}_{k} = \{y \in \mathbb{R}^{3} \, | \, y_{k} = 1\},
\quad
\overline{V}_{k} = \{y \in \mathbb{R}^{3} \, | \, y_{k} = -1\} 
\]
 for each local chart $U_{k}$ and $V_{k}$.
 We denote by $(x,\lambda)$ the value of $g^{\pm}_{k}(y)$ for any $k$.
 
\begin{figure}[ht]
\begin{center}
\includegraphics[width=8cm]{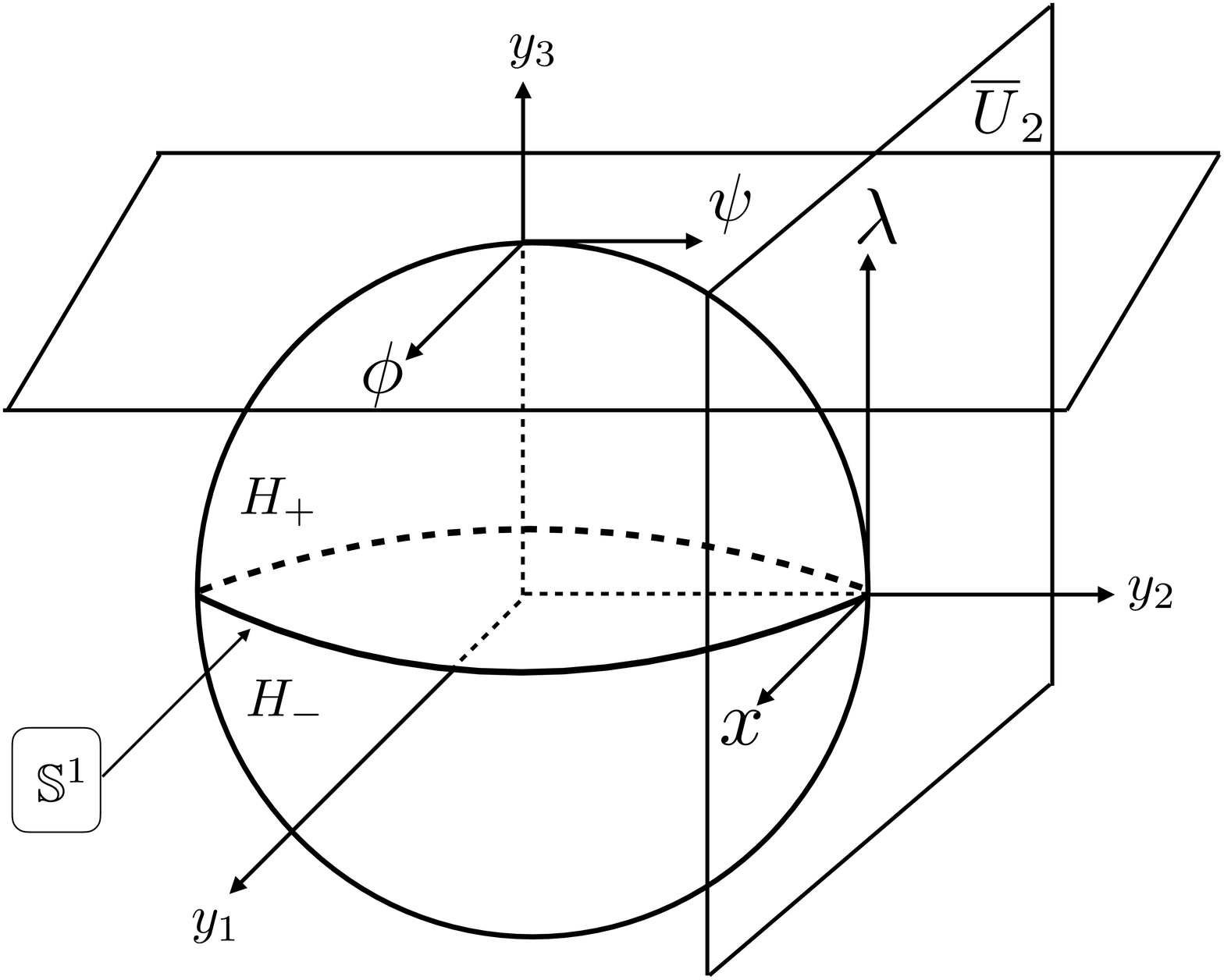}
\caption{Locations of the Poincar\'e sphere and chart $\overline{U}_{2}$.}
\label{fig:Poincare}
\end{center}
\end{figure}

For instance, it follows that
 \[ (g^{+}_{2} \circ f^{+})(\phi,\psi) = \left ( \dfrac{\phi}{\psi},\dfrac{1}{\psi}\right) = (x,\lambda), \]
therefore, we can obtain the dynamics on the local chart $\overline{U}_{2}$ by the change of variables $\phi = x/\lambda$ and $\psi = 1/\lambda$.
The locations of the Poincar\'e sphere, $(\phi,\psi)$-plane and $\overline{U}_{2}$ are expressed as Figure \ref{fig:Poincare}. 
Throughout this paper, we follow the notations used here for the Poincar\'e compactification.
 It is sufficient to consider the dynamics on $H_{+}\cup\mathbb{S}^{1}$, which is called Poincar\'e disk.
 
\section*{Appendix B: Two-dimensional ODEs characterizing the traveling wave solutions of \eqref{eq:pDNPE1-po}}
\label{appendix:ap2}
\setcounter{figure}{0}
\renewcommand{\thefigure}{B.\arabic{figure}}
In this paper, the classification of nonnegative weak traveling wave solutions of \eqref{eq:pDNPE1-po} is obtained via \eqref{eq:pDNPE1-1}.
As mentioned in Section \ref{sec:pDNPE1}, it is difficult to grasp the dynamics of the two-dimensional ODEs characterizing the traveling wave solutions of \eqref{eq:pDNPE1-po}.
In this section, we describe the reasons why it is difficult to understand the dynamics.

By introducing the traveling wave coordinate 
\[
\phi(\xi) =V(t, x),\quad \xi=x-ct, \quad c>0
\]
into \eqref{eq:pDNPE1-po} and setting $\psi(\xi)=d(\phi^{m})/d\xi$, we obtain the following equations:
\[
\begin{cases}
\phi'=m^{-1}\phi^{1-m}\psi,
\\
\psi'=-m^{-1}c\phi^{1-m}\psi
\end{cases}
\quad \left(\, '=\dfrac{d}{d\xi}\right)
\]
In order to be able to handle $\phi=0$ and to prevent the appearance of trajectories passing through $\phi=0$, we desingularize it by the time-scale desingularization: 
\[
\dfrac{ds}{d\xi}=\dfrac{1}{m} \phi^{-m}
\]
Then we have
\[
\begin{cases}
\phi' =\phi\psi,
\\
\psi'=-c\phi\psi
\end{cases}
\quad \left(\, '=\dfrac{d}{ds}\right)
\]
From a simple calculation similar to the one in this paper, it is expected that the entire $\phi$ and $\psi$ axes are equilibria for this system, and that a center manifold exists in the neighborhood of each point on the $\phi$ and $\psi$ axes.
However, obtaining an approximation of the center manifold at each point on both axes requires careful analysis and is an open problem.
That is, the results on the classification of non-negative weakly progressive wave solutions of \eqref{eq:pDNPE1-po} are not obtained by deriving the dynamics of the above two-dimensional ODEs characterizing the progressive wave solutions of \eqref{eq:pDNPE1-po}, but it is wiser to choose a method via \eqref{eq:pDNPE1-1} at present as in the present paper. 


\begin{thebibliography}{99}

\bibitem{AFJ} 
\'Alvarez, M.J., Ferragut, A., Jarque, X.: 
A survey on the blow up technique, 
Internat. J. Bifur. Chaos Appl. Sci. Engrg, {\bf21}, 3108--3118 (2011).


\bibitem{Anada} 
Anada, K., Ishiwata, T.: 
Blow-up rates of solutions of initial-boundary value problems for a quasi-linear parabolic equation, 
J. Differential Equations, {\bf262}, 181--271 (2017).

\bibitem{AIU22}
Anada, K., Ishiwata, T., Ushijima, T.:
Asymptotic expansions of traveling wave solutions for a quasilinear parabolic equation,
Jpn. J. Ind. Appl. Math, {\bf{39}}, 889--920 (2022).

\bibitem{ange91}
Angenent, S.,
On the formation of singularities in the curve shortening flow,
J. Diff. Geom., {\bf{33}}, 601--633 (1991).

\bibitem{Ange}
Angenent, S., Vel\'azquez, J.J.L.:
Asymptotic shape of cusp singularities in curve shortening,
Duke Math. J., {\bf{77}}, 71--110 (1995).

\bibitem{Aron}
Aronson, D.G.:
The porous medium equation, in ``Some Problems in Nonlinear Diffusion'' (A. Fasano and M. Primjcerio, Eds.), Lecture Notes in Math., Springer-Verlag, New York/Berlin, 1986.

\bibitem{MB}
Brunella, M.: 
Topological equivalence of a plane vector field with its principal part defined through Newton polyhedra,  
J. Differential Equations,  {\bf85}, 338--366 (1990).

\bibitem{carr} 
Carr, J.: 
Applications of centre manifold theory, 
Springer-Verlag, New York-Berlin (1981).


\bibitem{FAL} 
Dumortier, F., Llibre, J., Art\'es, C.J.: 
Qualitative theory of planar differential systems, 
Springer-Verlag, Berlin (2006).

\bibitem{QTW} 
Ichida, Y., Sakamoto, T.O.: 
Quasi traveling waves with quenching in a reaction-diffusion equation in the presence of negative powers nonlinearity, 
Proc. Japan Acad. Ser. A Math Sci. {\bf96}, 1--6 (2020).

\bibitem{DNLA} 
Ichida, Y., Matsue, K., Sakamoto, T.O.: 
A refined asymptotic behavior of traveling wave solutions for degenerate nonlinear parabolic equations, 
JSIAM Lett. {\bf{12}}, 65--68 (2020).

\bibitem{DNPE} 
Ichida, Y., Sakamoto, T.O.: 
Traveling wave solutions for degenerate nonlinear parabolic equations, 
J. Elliptic Parabol. Equ. {\bf{6}}, 795--832 (2020). 

\bibitem{cDNPE}
Ichida, Y: 
Classification of nonnegative traveling wave solutions for the 1D degenerate parabolic equations, 
Discrete Contin. Dyn. Syst., Ser. B, {\bf{28}}, no. 2, 1116--1132 (2023).

\bibitem{CK} 
Kuehn, C.: 
Multiple Time Scale Dynamics, 
Springer, Berlin (2015).


\bibitem{LPT12}
Lin, Y.C., Poon, C.C., Tsai, D.H., 
Contracting convex immersed closed plane curves with slow speed of curvature,
Transactions of AMS, {\bf{364}}, 5735--5763 (2012).

\bibitem{Matsue1} 
Matsue, K.: 
On blow-up solutions of differential equations with Poincar\'e-type compactificaions, 
SIAM J. Appl. Dyn. Syst., {\bf17}, 2249--2288 (2018). 

\bibitem{Matsue2} 
Matsue, K.: 
Geometric treatments and a common mechanism in finite-time singularities for autonomous ODEs, 
J. Differential Equations, {\bf267}, 7313--7368 (2019).


\bibitem{PaVa}
De Pablo, A., V\'azquez, J.L.:
Traveling waves and finite propagation in a reaction-diffusion equation,
J. Differ Equ., {\bf{93}}, 19--61 (1991).

\bibitem{Poon} 
Poon, C.C.: 
Blowup rate of solutions of a degenerate nonlinear parabolic equation, 
Discrete Contin. Dyn. Syst. Ser. B, {\bf24}, 5317--5336 (2019).

\bibitem{Wiggins}
Wiggins, S.: 
Introduction to Applied Nonlinear Dynamical Systems and Chaos, 
Springer-Verlag, New York, 2003.

\bibitem{Wink03}
Winkler, M.:
Blow-up of solutions to a degenerate parabolic equation not in divergence form,
J. Differential Equations, {\bf{192}}, 445--474 (2003).


\bibitem{XuYin}
Xu, T., Yin J.:
Traveling waves in degenerate diffusion equations,
Commun. Math. Res., {\bf{39}}, 36--53 (2023).

\end{thebibliography}
\end{document}